\documentclass{article}%
\usepackage{epsfig}
\usepackage{enumerate,graphicx}
\usepackage{amsfonts}
\usepackage{graphicx}
\usepackage{amsmath,amsthm,amssymb}
\usepackage{amssymb}%
\DeclareGraphicsExtensions{.pdf, .jpg, .png , .bmp, .ps}
\usepackage{tikz}
\providecommand{\U}[1]{\protect\rule{.1in}{.1in}}
\newtheorem{thm}{Theorem}[section]

\newtheorem{prop}{Proposition}[section]

\theoremstyle{remark}
\newtheorem{rem}[thm]{\bf Remark}

\def\N{{\mathbb N}}

\def\R{{\mathbb R}}

\def\E{{\mathbb E}}

\def\1{{\mathbbm{1}}}
\def\lip{{\mbox{Lip}\,}}

\def\1{{{\mbox{${\rm{1\negthinspace\negthinspace I}}$}}}}

\def\cov{\mathop{\rm Cov}\limits}
\newcommand{\eref}[1]{(\ref{#1})}

\newcommand\beq{\begin{equation}}
\newcommand\eeq{\end{equation}}

\begin{document}

\title{Deviation inequalities for contractive infinite memory processes}

\author{
Paul Doukhan\footnote{ AGM UMR CNRS 8088, CY Paris Cergy University, France.
Email: doukhan@cyu.fr}\footnotemark[1],   \ \   Xiequan Fan\footnote{School of Mathematics and Statistics, Northeastern University at Qinhuangdao, Qinhuangdao, China. E-mail: fanxiequan@hotmail.com}\footnotemark[2] 
}
\date{}

\maketitle

\abstract{In this paper, we introduce a class of processes that contains many natural examples.
The interesting feature of such type processes is based upon its infinite memory that allows it to record a quite ancient history.
Then, using the martingale decomposition method,  we establish some deviation and moment  inequalities for
 separately Lipschitz functions of such a process,
 under various moment conditions on some dominating random variables.  Our results generalize
 the Markov models of  Dedecker and Fan [Stochastic Process.\ Appl., 2015]  and a recent paper by Chazottes et al.\ [Ann.\ Appl.\ Probab., 2023]
  for the special case of a specific class of infinite memory models with discrete values.
   An application to stochastic gradient Langevin dynamic
 algorithm is also discussed.

\medskip


\medskip


\section{Introduction}\label{S1}
Concentration inequalities are essential tools for proving consistency and ensuring the validity of many statistical procedures;
let us cite for instance aggregation or selection model procedures as in \cite{A24} or prediction, see for example  \cite{AW12}.

In this paper, we consider a very general class of processes that includes many natural examples, such infinite memory models were introduced in \cite{IT69} for discrete state space models, and we aim at advocating the use of general state space models.  We first wanted to provide some comments to justify those models.
For each stationary and ergodic process $(X_n)_{n \in \mathbb{Z}}$, we denote ${\cal F}_{a}^b$ the $\sigma$-field generated by $X_n$ for $n\in (a,b)$. Typically, we may  write
$$
X_t=Z_t+\xi_t, \qquad Z_t=\E (X_t|{\cal F}_{-\infty}^{t-1}).
$$
Hence the process $\xi_t$ is such that
$\E ( \xi_t|{\cal F}_{-\infty}^{t-1})=0.$
As an  example of this situation, consider  stochastic volatility models $\xi_t$ with $\xi_t=\eta_tY_t$ with $Y_t$ an ${\cal F}_{-\infty}^{t-1}-$measurable random variable and $(\eta_t)_t$ an independent and identically distributed (i.i.d.) centered sequence to see that the previous decomposition may not be rich enough to model the memory of a process.
 Such models were in fact generalised in \cite{DW08}  and led to infinite memory models.
The interesting feature of such type processes relies on its infinite memory that allows it to record a quite ancient history.
Notice that infinite memories models may  approximate by finite memory models since:
$$
\E (X_t|{\cal F}_{-\infty}^{t-1})=\lim_{d\to\infty}\E (X_t|{\cal F}_{t-d}^{t-1})$$
which is also a natural justification of Markov modelling.
A fruitful approach to such models is in  \cite{BW99} who consider Markov chains with a memory depending on the past of the process; those authors  also prove  that infinite memory models are natural models of time series.

Concerning now the structure of the paper, using the martingale decomposition method,  we establish some deviation and moment  inequalities for
 separately Lipschitz functions of such a process,
 under various moment conditions on some dominating random variables.
 Several papers were previously devoted to $d$-th order Markov models  (cf. \cite{DF15,DDF19})  and a recent paper \cite{CGT23}  considers the special case of a specific class of infinite memory models with discrete values.
As it was mentioned in these paper, many applications of such models are natural. 

\section{Iterated random functions with infinite memory}
In this section, we introduce a model with infinite memory. Some explicit examples of such type models  are also presented.
$d$-order Markov models may be simply rephrased as solutions of a recursion
$$
X_n=F_n(X_{n-1},\ldots, X_{n-d},\xi_n)
$$
under very weak assumptions on the state space.
The problem to address very large values of the memory lead to consider $d=\infty$, or even $d$ arbitrary large.

\subsection{An infinite memory process}
Let $(\Omega, {\mathcal A}, {\mathbb P})$ be a probability space.
Let $({\mathcal X}, d)$ and  $({\mathcal Y}, \delta)$ be two complete separable metric spaces.
Let $(\varepsilon_i)_{i \geq 1}$ be a sequence of independent (not necessary identically distributed)  ${\mathcal Y}$-valued
random variables. Let $(X_i)_{i\leq 0}$ be a ${\mathcal X}^{\mathbb{N}}$-valued random variable independent of $(\varepsilon_i)_{i \geq 1}$. We consider the  (non-Markov) infinite memory processes  $(X_i)_{i \geq 1}$
such that
\begin{equation}\label{Mchain}
X_n=F_n ( (X_{n-i})_{i\geq 1};   \varepsilon_n), \quad \text{ $n\geq 1$},
\end{equation}
where $F_n: {\mathcal X}^{\mathbb{N}} \times {\mathcal Y} \rightarrow {\mathcal X}$
is such that
\begin{equation}\label{contract}
{\mathbb E}\big [ d\big(F_n ( (x_{n-i})_{i\geq 1}   ;   \varepsilon_n), F_n ((x'_{n-i})_{i\geq 1}  ; \varepsilon_n)\big) \big] \leq
 \sum_{i= 1}^{ \infty} a_{i} \,   d(x_{n-i}, x_{n-i}')
\end{equation}
for some nonnegative numbers $a_{i}, i\geq 1,$ such that $ \sum_{i= 1}^{ \infty} a_{i} < 1$.
{\color{black}When $F_n\equiv F$, such contractive model  \eqref{Mchain} is introduced in   \cite{DW08}   and additionally to \eqref{contract} an existence and uniqueness condition additionally needs: for some $A\in{\mathcal Y}$,
$$
{\mathbb E}\big [ d\big(F_n  ((x_{n-i})_{i\geq 1}; \varepsilon_n),   A\big) \big]<\infty.  $$}


\begin{rem} Let us give some typical models of type (\ref{Mchain}).
\begin{enumerate}
\item  \emph{Markov models.}  Dedecker and Fan \cite{DF15} considered the following Markov model
\begin{equation}\label{fdsf32}
X_n=F(  X_{n-1} ;   \varepsilon_n), \quad \text{ $n\geq 1$},
\end{equation}
where $F : {\mathcal X}  \times {\mathcal Y} \rightarrow {\mathcal X}$
is such that
\begin{equation} \label{contract00}
{\mathbb E}\big [ d\big(F ( x_{n-1}   ;   \varepsilon_n), F ( x'_{n-1}  ; \varepsilon_n)\big) \big] \leq
  \rho \, d(x_{n-1}, x_{n-1}')
\end{equation}
for some $\rho \in [0, 1).$  It is known that there are a lot of models satisfying condition \eqref{contract00}, see   \cite{DF15}. Clearly,  condition (\ref{contract00}) implies (\ref{contract}) with $F_n(   (x_{n-i})_{i\geq 1} ;   \varepsilon_n)=F(x_{n-1} ; \varepsilon_n),$ $a_1= \rho$ and $a_i=0, i\geq 2$. Thus, our model is an extension of the model (\ref{fdsf32}).

\item\emph{A Non-Markov model.}
Consider thus  the following model:
\begin{equation}
X_n=F_n  (  X_{n-1}, \ldots, X_1 ;   \varepsilon_n), \quad \text{ $n\geq 2$},
\end{equation}
where $F_n : {\mathcal X}^{n-1}  \times {\mathcal Y} \rightarrow {\mathcal X}$
is such that
\begin{equation} \label{contract02}
{\mathbb E}\big [ d\big(F_n ( x_{n-1},\ldots,x_1   ;   \varepsilon_n), F_n ( x'_{n-1},\ldots,x'_1  ; \varepsilon_n)\big) \big] \leq
 \sum_{i= 1}^{ n-1} a_{i}\,  d(x_{n-i}, x_{n-i}')
\end{equation}
for some $a_i \geq 0, i \geq 1,$ and $ \sum_{i= 1}^{ \infty} a_{i} < 1$. Clearly, condition (\ref{contract02}) implies (\ref{contract}) with $$F_n(   (x_{n-i})_{i\geq 1} ;   \varepsilon_n)=F_n(x_{n-1},...,x_1; \varepsilon_n).$$
The time non homogeneous case is in  \cite{DDF19} and $F_n(x_{n-1},...,x_1; \varepsilon_n)=G_n(x_{n-1}; \varepsilon_n)$ with $G_n=G$ in
  \cite{DF15} and $G_n$ varies with $n$ in   \cite{DDF19}.
 Quote that this is a bit different model since $F_n$ is defined here on $ {\mathcal X}^{n-1} $ contrary to \eqref{Mchain} where $F_n$ was defined over the bigger product set
 $ {\mathcal X}^\N$; this latter model may be defined recursively from time $n=0$ and thus corresponds better to practitionner's intuitions.

\item\emph{$p$-Markov models.}   Model (\ref{contract}) also includes the following $p$-Markov model
\begin{equation}
X_n=F(  X_{n-1},\ldots, X_{n-p};   \varepsilon_n), \quad \text{ $n\geq 1$},
\end{equation}
where $F : {\mathcal X}^p  \times {\mathcal Y} \rightarrow {\mathcal X}$
is such that
\begin{equation} \label{contract01}
{\mathbb E}\big [ d\big(F (x_{n-1},\ldots, x_{n-p} ;   \varepsilon_n), F ( (x'_{n-1},\ldots, x'_{n-p}  ; \varepsilon_n)\big) \big] \leq
 \sum_{i= 1}^{ p} a_{i}  d(x_{n-i}, x_{n-i}')
\end{equation}
for some $a_i \geq 0, i \geq 1,$ and $ \sum_{i= 1}^{ p} a_{i} < 1$. Clearly, condition (\ref{contract01}) implies (\ref{contract}) with $$F_n(   (x_{n-i})_{i\geq 1} ;   \varepsilon_n)=F (x_{n-1},...,x_{x-p}; \varepsilon_n).$$
\end{enumerate}
\end{rem}

The main attraction of the model (\ref{Mchain}) lays on the infinite memory. Thanks to the  memory effect (\ref{contract}),
model (\ref{Mchain}) is allowed to record quite  a lot of  history. Therefore, it provides a useful tool to model
 data that exhibit infinite memories.

\subsection{Examples}\label{secexample}
In this subsection, we give a non exhaustive list of models satisfying condition \eqref{contract}.
\begin{enumerate}

\item \emph{ARCH-type  models}. For the model
$$
X_n= \sqrt{\sum_{i=1}^\infty a_i^2X_{n-i}^2   +b^2} \cdot \varepsilon_n,
$$
set $$F_n(x;y)= F (x;y)=  \sqrt{\sum_{i=1}^\infty a_i^2 x_{i}^2   +b^2} \cdot y \ \ \ \ \textrm{and} \ \ \ \ d(x,x')= |x-x'|.$$
Then, it is easy to verify that
\begin{equation}
{\mathbb E}\big [ d\big(F  ( (x_{n-j})_{j\geq 1}   ;   \varepsilon_n), F  ((x'_{n-j})_{j\geq 1}  ; \varepsilon_n)\big) \big] \leq
 \sum_{j= 1}^{ \infty} |a_j| \,  \mathbb{E} |\varepsilon_n| \, d(x_{n-j}, x_{n-j}').
\end{equation}
For this model, contraction  \eqref{contract} is satisfied   provided that  $  \sum_{j= 1}^{ \infty} |a_j|  \sup_n {\mathbb E}  |\varepsilon_n| < 1.$

\item
\emph{GLM type models}.   Assume that $(\xi_n)_n$ is an i.i.d.\ sequence of unit Poisson processes
$$X_n=\xi_n(\lambda_n),\qquad
\lambda_n=g(X_{n-t_1}, X_{n-t_2},..., ,X_{n-t_p})$$
with a contractive function $g$ such that
$$|g(y_1,y_2,\ldots, y_p)-g(x_1,x_2,\ldots, x_p)|_r\le \sum_{j=1}^p a_j|y_j-x_j|, \qquad \sum_{j=1}^p a_j<1.$$   

\item \emph{Memory one/infinite models.}  Assume that $(J_n)_{n\geq 1}$ is an i.i.d.\ sequence. Let
$$
X_n=a_{J_n} X_{n-J_n}+\xi_n.
$$
In this case $\sum_{i=1}^{\infty} |a_i| {\mathbb P}(J=i)<1$ is a stationarity condition.
Quote that $a_i\equiv1$ is not a possible value, but it is enough that one of the coefficients  be $|a_{i_0}|<1$ and it is also possible to consider some explosive regimes $|a_i|>1$ in case they are compensated by contractive ones, for instance, ${\mathbb P}(J=i) < (1 - \sum_{j\neq i }^{\infty} |a_j| ) /|a_i|$.
Thus, contraction (\ref{contract}) holds.

\item Consider  this  extension of the elephant walk:
$$X_n=A_n X_{n-Z_n}+\zeta_n,$$
with $\xi_n=(A_n, Z_n)$ an i.i.d.\ sequence such that $\E |A_n|^p<1$ and  $\E |\zeta_n|^p<1$ and  $Z_n \in \{1, 2, 3, \ldots\} $, then it fits our condition (\ref{contract}).

\item \emph{Generalized elephant random walks.}
 Let
$(\xi_i)_{i\geq 1}$ be a sequence of i.i.d.\ random variables.
The random walk starts at the origin at time
zero, $S_0=0$.   At time $n = 1$, the walker moves to $\xi_1.$   Hence, the position of the walker
at time $n = 1$ is given by $S_1 =X_1$ with $ X_1=\xi_1$. Afterwards,
at any time $n\geq 2$, we choose the random variable $X_n$ with probability $t, t \in [0, 1],$ or choose  at random an integer $k$ among the previous times $1, 2, \ldots, n-1$ with equal probability
$ (1-t)/(n-1)$. Define
\begin{displaymath}
X_{n} = \left\{ \begin{array}{ll}
\xi_n & \textrm{  with probability $t $ }\\
 X_k & \textrm{  with probability $(1-t)p $  }\\
-X_k & \textrm{  with probability $(1-t)(1-p) $},
\end{array} \right.
\end{displaymath}
where the parameter  $p \in [0, 1] $ is the memory of the ERW. Then, the position of the ERW is
given by $$S_{n} = S_{n-1} + X_{n }.$$

In order to understand well how the elephant moves, it is straightforward to see that for any
time $n \geq 2,$ $$X_{n } = \gamma_n X_{\eta_n} + (1-|\gamma_n|)\xi_n$$ where $ \gamma_n, \eta_n$ and $\xi_n$ are   independent discrete random variables,
with the distribution
\begin{displaymath}
\gamma_{n} = \left\{ \begin{array}{ll}
  0& \textrm{with probability $t  $  }\\
 1 &  \textrm{with probability $(1-t)p$}\\
-1 & \textrm{with probability $(1-t)(1-p)$}
\end{array} \right.
\end{displaymath}
 and $\eta_n$ is uniformly  distributed over the integers
$\{ 1,  \ldots , n-1  \}$ such that $\mathbb{P}(\eta_n=k )= (1-t)/(n-1) $. Moreover, $\{\gamma_n\}_{n\geq 1}$ is independent of  $\{X_n\}_{n\geq 1}.$
Clearly, when $t=1,$ $S_n$ reduces to the elephant random walk, see \cite{ST04}.
When $p=1,$ $S_n$ reduces to the step-reinforced random walk, see \cite{B21}.
It is easy to see that
\begin{equation}
{\mathbb E}  | \gamma_n x_{\eta_n}-\gamma_n x'_{\eta_n}|   \leq
 \sum_{k= 1}^{ n-1} \frac{1-t}{n-1}   |x_{n -k} - x'_{n -k}|
\end{equation}
and
$$\sum_{k= 1}^{ n-1} \frac{1-t}{n-1}   =  1-t < 1,$$
provided that $t>0.$
Therefore, condition (\ref{contract}) is satisfied.

\item \emph{Random memory AR-models.}
Assume that coefficients $a_1,a_2,\ldots$ satisfy $ \sum_{i=1}^\infty |a_i|<1$ and that $\varepsilon=(J,\xi)\in {\mathbb N}\times {\mathbb R}$ then from an i.i.d.\ sequence. With this distribution, it is easy to define
$$
X_n=\sum_{i=1}^{J_n} a_iX_{n-i}+\xi_n.
$$
Here with $
F(x;\varepsilon)=\sum_{i=1}^{J} a_ix_{i}+\xi
$, contraction (\ref{contract}) holds.

\item \emph{Mean fields memory models.} Assume that
$$
X_n=r\bigg(\sum_{i=1}^{\infty} a_iX_{n-i}\bigg)+\varepsilon_n.
$$
In case the coefficients $a_i$ are known then the model is a simple regression model
$$ X_n=r(Y_n)+\varepsilon_n,
\qquad
Y_n=\sum_{i=1}^\infty a_i X_{n-i}.$$
Then
$$
F(x;z)=r\Big(\sum_{i=1}^\infty a_i x_{i} \Big)+z.
$$  It is easy to see that
\begin{eqnarray*}
{\mathbb E}|F(x';y)-F(x;y)|\le \, |x'-x | \, \lip r \ \le\    \lip r  \sum_{i=1}^\infty |a_{i}|   |x'_{i}-x_{i}|,
\end{eqnarray*}
where we set $x=\sum_{i=1}^\infty a_i x_{i}$ and  $x'=\sum_{i=1}^\infty a_i x'_{i}$.
Then, contraction (\ref{contract}) holds in case $  \sum_{i=1}^\infty |a_{i}| \in [0,  1/ \lip r   )$.

%

%

\item Consider two bounded functions
$$
\varphi:{\mathbb R}\to{\mathbb R} \ \ {\color{black}\textrm{and} \ \  \psi:{\mathbb R}\to{\mathbb R^+}}
$$
with $\|\varphi\|_\infty,\|\psi\|_\infty<\infty$,
then a neural based model writes
\begin{eqnarray}
\label{model}
X_t&=&\sum_{j=1}^{Z_t}a_j\varphi(X_{t-j})+\xi_t,\\
\nonumber\\
Z_t&=&P (\lambda_t),\\
\lambda_t&=&\psi(X_{t-1}),
\end{eqnarray}
{\color{black} where $P(\lambda)$ is a Poisson point process.}
The main attraction of the model  is that the geometric memory effect makes it recording about a quite ancient history.
First quote that $X_t=F(X_{t-1}, X_{t-2},\ldots;\zeta_t)$ with $\zeta_t=(\xi_t,P(\lambda_t) )$
and $$
F(x;\zeta)=\sum_{j=1}^{Z}a_j\varphi(x_{j})+\xi,\
Z=P(\lambda),\
\
\lambda=\psi(x_{1}),
$$
where $x=(x_1, x_2, ... ).$
Hence, if ${\mathbb E}|\xi|<\infty,$ then
$$
{\mathbb E} |F(0,\zeta)|\le  \|\varphi\|_\infty \sum_{j=1}^{\infty} |a_j| +{\mathbb E}|\xi|<\infty.
$$
Now if $\psi(x'_1)\ge \psi(x_1)$ the monotonicity of $\lambda\mapsto P(\lambda)$ implies, with $Z=P(\psi(x_{1}))$ and $Z'=P(\psi(x'_{1}))$, that
\begin{eqnarray*}
F(x';\zeta)-F(x;\zeta)&=&\sum_{j=1}^{Z}a_j(\varphi(x'_{j})-\varphi(x_{j}))+\sum_{j=Z+1}^{Z'}a_j\varphi(x'_{j}).
\end{eqnarray*}
From the last inequality, we deduce that
\begin{eqnarray*}
{\mathbb E}|F(x';\zeta)-F(x;\zeta)|&\le&\sum_{j=1}^{\infty}|a_j| |\varphi(x'_{j})-\varphi(x_{j})|+\sup_j|a_j|\|\varphi\|_\infty{\mathbb E}|Z'-Z|
\\
&\le&  \lip \varphi \sum_{j=1}^{\infty} |a_j| |x'_{j}-x_{j}|+\sup_j|a_j| \,\|\varphi\|_\infty\lip \psi  |x'_1-x_1|
\\
&\le& ( |a_1| \lip \varphi + \sup_j|a_j| \,\|\varphi\|_\infty\lip \psi) |x'_{1}-x_{1}|   + \sum_{j=2}^{\infty}|a_j| |x'_{j}-x_{j}|.
\end{eqnarray*}
 Thus, condition (\ref{contract}) is satisfied, provided that
  $$  |a_1|  \lip \varphi + \sup_j|a_j| \,\|\varphi\|_\infty\lip \psi    + \sum_{j=2}^{\infty}|a_j| < 1.$$

\end{enumerate}

\section{Separately Lipschitz functions}
\setcounter{equation}{0}
For each integer $n\ge1$, let $f: {\mathcal X}^n \to {\mathbb R}$ be a separately Lipschitz function, such that
\begin{equation} \label{codiMD}
|f(x_1, x_2, \ldots, x_n)-f(x'_1, x'_2, \ldots, x'_n)| \leq   d(x_1,x'_1)+ d(x_2,x'_2)+ \cdots +     d(x_n, x'_n) \, .
\end{equation}
The natural filtration of the chain is defined as  ${\mathcal F}_0=\{\emptyset, \Omega \}$
and, for $k \in {\mathbb{N}}^{*}$,
${\mathcal F}_k= \sigma(X_1, X_2,  \ldots, X_k)$.
Define
\begin{equation}\label{gk}
g_k(X_1, X_2, \ldots , X_k):= {\mathbb E}[f(X_1, X_2, \ldots, X_n)|{\mathcal F}_k]\, .
\end{equation}
Denote
\begin{equation}\label{dk}
d_k=g_k(X_1, X_2, \ldots, X_k)-g_{k-1}(X_1, X_2, \ldots, X_{k-1})
\end{equation}
and
\begin{equation}\label{Sn}
S_n =f(X_1, X_2,\ldots , X_n) -{\mathbb E}[f(X_1, X_2, \ldots , X_n)]\,.
\end{equation}
It is easy to see that $(d_k, {\mathcal F}_k)_{k=1,2,\ldots,n}$ is a finite sequence of martingale differences, For $k \in [1, n-1]$, let
$$
S_k:=d_1+d_2+\cdots + d_k .
$$
Thenl $S_n =d_1+d_2+\cdots + d_n   .
$.
By the definition of  $d_k$'s,  this is easy to see that   $(S_k, {\mathcal F}_k)_{k=1,2,\ldots,n}$ is a martingale.

The following proposition gives  some interesting properties of
the functions $g_k$ and of the martingale differences $d_k$.

\begin{prop}\label{McD}
Let $k \in {\mathbb N},$ and let $(X_i)_{i \geq 1}$ be a   chain satisfying \eref{Mchain} for some
functions $F_n$ satisfying \eref{contract}. Let $g_k$ and $d_k$ be defined
by
\eref{gk} and \eref{dk} respectively.
\begin{enumerate}
\item
The function $g_k$ is separately Lipschitz and satisfies
$$
\Big|g_{k}(x_1, x_2,  \ldots, x_{k})-g_{k}(x'_1, x'_2, \ldots, x'_{k}) \Big| \leq  \sum_{l=1}^{k} a_{ n-k}(n-l) \, d(x_{ l}, x'_{ l})   ,
$$
where $$a_0(0)=1,  \, \  a_0(i)=1, \, \ a_1(i)=1+a_i,\ \ \ \ \ \ \ \ \ \ \ \ \ \ \ \ \ \ \ \ \ \ \ \ \ \ \ \ \ \ \ \ \ \ \ \ $$
$$   a_{k+1}(i)=a_{k}(i)+ a_{k}(k) a_{i-k}, \ \ \  k \in [1, \, n-1],\  \textrm{and} \ i \in [k+1, n-1].$$
In particular,  we have
\begin{eqnarray}\label{akdefd}
a_0(0)=1,\ \ \ a_1(1)=1+a_1,\ \ \    a_k(k)= 1+a_k + \sum_{l=1}^{k-1} a_l(l)a_{k-l}.
\end{eqnarray}
\item Let $\widetilde{P} $ be the distribution of the random vector $(X_i )_{i\leq 0}$ and $P_{\varepsilon_k}$ be the
distribution of the $\varepsilon_k$'s. Denote by $\mathbf{X}_k=(X_k, X_{k-1},\ldots,X_1).$
Let  $H_{k,\varepsilon_k}$ be  defined by for $k=1,$
$$
H_{1, X_1}( X_1) =  \int  d(X_1,
 y)P_{X_1}(dy) \,
$$
and for any $k \in [2,n]$,
\begin{eqnarray*}
  H_{k, \varepsilon_k}(\mathbf{X}_{k-1},  \varepsilon_k)      =\! \int\!\!\!\!\int \!\! d(F_k(\mathbf{X}_{k-1},x_0,,\ldots; \varepsilon_k),
 F_k(\mathbf{X}_{k-1},x_0,\ldots; y))P_{\varepsilon_k}(dy)\widetilde{P} ( dx_0,\ldots).
\end{eqnarray*}
Then, the martingale difference $d_k$ satisfies
$$
   |d_1|\leq a_{n-1}(n-1)  H_{1,X_1}( X_1 )   \ \ and   \ \ |d_k|\leq a_{n-k}(n-k)  H_{k,\varepsilon_k}(\mathbf{X}_{k-1},  \varepsilon_k), \   k \in [2,n].
$$
\item
Assume moreover that $F_n$ satisfies
\begin{equation}\label{c2}
d(F_n((x_i)_{i\leq k};y), F_n((x_i)_{i\leq k};y')) \leq   \delta(y,y'),
\end{equation}
 and
let
$$
G_{X_1}(y)= \int  d(y,
 y')P_{X_1}(dy')\ \  \textrm{and} \ \   G_{\varepsilon_k}(y)= \int  \delta(y,y')P_{\varepsilon_k}(dy'), \  k \in [2,n].
$$
Then, for any $k \in [2,n]$, $H_{k,\varepsilon_k}(\mathbf{X}_{k-1},  y) \leq G_{\varepsilon_k}(y)$,
and consequently,
$$
 |d_1|\leq a_{n-1}(n-1) G_{X_1}(X_1)\  \ and \ \ |d_k|\leq a_{n-k}(n-k) G_{\varepsilon_k}(\varepsilon_k), \  k \in [2,n].
$$

\item
Assume that $F_n$ satisfies (\ref{c2}),  and that
$(X_{i})_{i\leq 0}$ are deterministic.  Then, for any $k \in [1,n]$,
$$
 |d_k|\leq a_{n-k}(n-k) G_{\varepsilon_k}(\varepsilon_k)\, .
$$
\end{enumerate}
\end{prop}

\begin{rem}\label{re2.1}
 Let us give some comments on Proposition \ref{McD}.
\begin{enumerate}
\item We first comment equations  (\ref{akdefd}).  Without loss of generality, we may assume that $a_k$ is decreasing with respect to $k.$ This is always true since we may change the order of the sequence $(x_i)_{i\geq 1}$ in (\ref{contract}). It is easy to see that
\begin{eqnarray}
   a_{k+1}(k+1)\ \geq \ 1  + \sum_{l=1}^{k-1 } a_l(l)a_{k +1-l }+  a_k(k)a_{1 }\  \geq \ 1 + \sum_{l=1}^{k-1 } a_l(l)a_{k -l }+a_{k }\ = \ a_k(k), \nonumber
\end{eqnarray}
which implies that $a_k(k)$ is increasing in $k$. Thus we have
$$a_k(k) \ \leq \ 1+ a_k + a_k(k) \sum_{i=1}^{k-1} a_{k-i} \ \leq \ 1+ a_k + a_k(k) \sum_{i=1}^{\infty} a_{i}.$$
From the last line, we get, for all $k\geq 1,$
$$ a_{k}(k) \leq  \frac{1+   \max_{ i }a_i  }{1- \sum_{i= 1}^{ \infty} a_{i}}.$$
Recall that thanks to  equation (\ref{contract})  we have  $\sum_{i= 1}^{ \infty} a_{i} < 1$. Thus, the last line implies that $(a_{k}(k))_{k\geq 1}$ is uniformly bounded with respect to $k.$

\item
It is easy to see that all examples in subsection \ref{secexample} satisfy the point 2 of  Proposition \ref{McD}.
It is also easy to see that the example 7 in Subsection \ref{secexample} and stochastic gradient Langevin dynamic in Subsection \ref{SGLD}
satisfy condition \eqref{c2}.

\item Let us comment on the point 4 of Proposition \ref{McD}. The fact that for each integer $k$,  the martingale difference $d_k$ is bounded by the random variable
$ a_{n-k}(n-k)G_{\varepsilon_k}(\varepsilon_k)$
which is {\it independent  of }${\mathcal F}_{k-1}$ plays a crucial role.
With this insight, we find that for any positive and increasing function $f$, it holds
$\mathbb{E}[f(d_k) | \mathcal{F}_{k-1} ] \leq \mathbb{E}[f(a_{n-k}(n-k)G_{\varepsilon_k}(\varepsilon_k))].$
 This explains why we  obtain deviations inequalities for $S_n$ under some conditions on the
distribution of $G_{\varepsilon_k}(\varepsilon_k).$

\item If $a_1=\rho$ and $a_i=0$ for any $i\geq 2,$ by (\ref{akdefd}), then we have
\begin{eqnarray*}
a_0(0)=1 \ \ \    \textrm{and} \ \ \   a_k(k)= 1+ \rho + ...+ \rho^k, \ \ k\geq 1.
\end{eqnarray*}

\end{enumerate}
\end{rem}

\section{Deviation  inequalities for the functional $S_n$} \label{deviationiq}
\setcounter{equation}{0}

 In this section, we assume that $(X_{i})_{i\leq 0}$ are \emph{deterministic}. We present some deviation inequalities  for the  functional $S_n$,
  with $X_1, X_2, ...,X_n$ satisfying the assumptions (\ref{Mchain}) and  (\ref{contract}).
Thanks to Proposition  \ref{McD}, the proofs of these inequalities are close to
that of \cite{DF15,DDF19}.
We  present the proofs of the propositions of this section   in Appendix.

Let us note that the deviations
inequalities   of this section are given for ${\mathbb P}\big(\pm S_n >x\big)$,
but we shall only prove them for $S_n$. The proofs of the
deviation inequalities for $-S_n$ are exactly the same,  because  the upper bounds of  points 2 and 3 of Proposition \ref{McD} hold both for $d_k$ and $-d_k$.

When $(X_{i})_{i\leq 0}$ are random variables,
thanks to point 3  of Proposition \ref{McD}, the deviation inequalities in this section
hold also but  with $G_{\varepsilon_1}(\varepsilon_1 )$ replacing by $G_{X_1}(X_1 )$.

\subsection{Bernstein type bound}\label{Bersec}
The well-known Bernstein inequality gives a tight Gaussian type bound on tail probabilities for sums of independent random variables.
Under the conditional Bernstein inequality, van de Geer \cite{V95}  and de la Pe\~{n}a \cite{D99} have established
the generalizations of Bernstein inequality for martingales. Now, applying  Proposition \ref{McD}, we have the following Bernstein type inequality.
\begin{prop}\label{pr01}
Assume condition (\ref{c2})  and that there exist some constants  $M>0$ and $V_k\geq 0$   such that, for any $k \in [1, n]$ and any $l\geq 2$,
\begin{equation}\label{BernsteinC}
 {\mathbb E} \Big[  \big(  G_{\varepsilon_k}(\varepsilon_k )\big)^l\Big] \leq
 \frac {l!}{2} V_k M^{l-2} \, .
\end{equation}
Let $$V=  \sum_{k=1}^n \big(   a_{n-k}(n-k) \big)^2 V_k \ \ \ \ \ \textrm{and}\ \ \ \ \ \delta=M a_{n-1}(n-1).$$
Then, for any $t \in [0, \delta^{-1})$,
\begin{equation} \label{maindfs}
  \mathbb{E} e^ {\pm  tS_n}  \leq \exp \left (\frac{t^2 V }{2 (1- t\,\delta )} \right )\, .
\end{equation}
Consequently, for any $x> 0$,
\begin{eqnarray*}
  {\mathbb P}\Big(\pm  S_n \geq  x\Big)
  &\leq&  \exp \left(\frac{-x^2}{V(1+\sqrt{1+2x \delta/V})+x \delta }\right)\,   \\
 &\leq&  \exp \left(\frac{-x^2}{2 \,(V +x  \delta ) }\right)\, .
\end{eqnarray*}
\end{prop}


Without condition (\ref{c2}), we have the following result under a sub-Gaussian  type condition.
\begin{prop} \label{BerProp}
Assume that there exists a positive constant $\epsilon$   such that,
for any   $k\geq 1$ and any $l\geq  2,$
\begin{equation}\label{Bercond02}
\mathbb{E}\left[ \big(H_{k,\varepsilon_k}(X_{k-1},\ldots,X_1, \varepsilon_k) \big)^l  \right] \leq \frac{1}{2} \frac{  l! \,\epsilon ^{l-2}  }{(l-1)^{l/2}}  \mathbb{E}\left[  \big(H_{k,\varepsilon_k}(X_{k-1},\ldots,X_1, \varepsilon_k) \big)^2 \right] .
\end{equation}
 Then,  for any   $ x>0$,
\begin{align}
  \mathbb{P}\Big( \pm S_n \geq x V_n  \Big)
   &\leq  \exp\bigg( \frac{- x ^2}{  1+ \sqrt{ 1+ 2 x \epsilon  a_{n-1}(n-1)
  /\sigma_n   }  + x   \epsilon a_{n-1}(n-1)/\sigma_n   } \bigg) \ \label{jnsk01} \\
 &\leq   \exp\bigg( \frac{-  x ^2}{  2 \big( 1+ x   \epsilon a_{n-1}(n-1)  / \sigma_n   \big)} \bigg),\label{jnsk02}
\end{align}
where $$V_n^2 =  \sum_{k=1}^{n} (a_{n-k}(n-k))^2 \mathbb{E}\big[ \big( H_{k,\varepsilon_k}(X_{k-1},\ldots,X_1, \varepsilon_k) \big)^2 \big] \ \ \ \  and \ \  \ \ \sigma_n^2=  \frac1n V_n^2.$$
\end{prop}


\begin{rem} Let us give some comments on Proposition \ref{BerProp}.
\begin{enumerate}
\item Condition \eqref{Bercond02} is in fact a sub-Gaussian   condition. Using Taylor's expansion,  on can check that it is satisfied provided that
$$   \inf_{k} {\mathbb E}\big[ \big( H_{k,\varepsilon_k}(X_{k-1},\ldots,X_1, \varepsilon_k) \big)^2 \big]    >0  $$
and
$$  \sup_{k}  {\mathbb E}\Big[ \exp  \Big( c \, \big( H_{k,\varepsilon_k}(X_{k-1},\ldots,X_1, \varepsilon_k)\big)^2 \Big) \Big]  < \infty, $$
 where $c$ is a  positive constant not depending on $k$.

\item Assume that
$$ 0< \liminf_{n\rightarrow \infty}\sigma_n    \leq \limsup_{n\rightarrow \infty}\sigma_n < \infty.$$
Then $V_n$ is of order $\sqrt{n}$. By Remark \ref{re2.1}, $a_n(n)$ is uniformly bounded with respect to $n$.
Therefore, by (\ref{jnsk02}), we find that
\begin{eqnarray}
\ln \mathbb{P}\left( \pm S_n \geq  n  \right)   \leq  - C \sqrt{n}
\end{eqnarray}
for some positive constant $C$ not depending on $n$.
\end{enumerate}
\end{rem}

\subsection{Cram\'{e}r type bound}
When $G_{\varepsilon_k}(\varepsilon_k)$ satisfies the Cram\'{e}r condition,  we obtain the following inequality
similar to that of  \cite{LW09} for martingales under the conditional Cram\'{e}r condition.
\begin{prop}\label{cram}
Assume condition (\ref{c2}) and that there exist some constants $t_0>0 $ and $K_k\geq 1$  such that
\begin{equation}\label{laplace}
 {\mathbb E} \Big[ \exp \Big( t_0   G_{\varepsilon_k}(\varepsilon_k)\Big)\Big] \leq K_k,\ \ \   \ k\in [1, n]
 \, .
\end{equation}
Let $$K=\frac{2}{e^{2}}   \sum_{k=1}^{n} \Big(\frac {a_{n-k}(n-k)}{a_{n-1}(n-1)}\Big)^2 K_k   \ \ \ \ \ and \ \ \ \ \ \delta=\frac{t_0}{a_{n-1}(n-1)} .$$
Then, for any $t \in [0, \delta )$,
$$
  \mathbb{E}\,[e^{\pm  tS_n}]\leq   \exp \left( \frac{t^2 K  \delta^{-2} }{1-t \delta^{-1}}   \right) \, .
$$
Consequently, for any $x> 0$,
\begin{eqnarray}
  {\mathbb P}\big(\pm  S_n \geq  x\big) &\leq&  \exp \left( \frac{-(x\delta)^2}{2K (1+\sqrt{1+ x \delta / K })+x \delta  }\right)\,   \label{Berie3}  \\
 &\leq&  \exp \left( \frac{-(x\delta)^2}{ 4K+ 2 x  \delta  }\right)\,\, .\label{Berie4}
\end{eqnarray}
\end{prop}

Recall that, by point 1 in Remark \ref{re2.1}, $a_n(n)$ is uniformly bounded with respect to $n$.
Assume that
$$\sup_k {\mathbb E} \Big[ \exp \Big( t_0   G_{\varepsilon_k}(\varepsilon_k)\Big)\Big] < \infty.$$
Then $K$ is of order $n$ as $n\rightarrow \infty$.
Therefore, by Proposition \ref{cram}, it is easy to see that
\begin{eqnarray}
\ln \mathbb{P}\left( \pm S_n \geq  n  \right)   \leq  - C n
\end{eqnarray}
for some positive constant $C$ not depending on $n$.

\subsection{Semi-exponential bounds}
When  $G_{\varepsilon_k}(\varepsilon_k)$
has semi-exponential moment, we have the  the following proposition. This proposition can be compared to the
corresponding results in \cite{Bor} for   sums of independent random variables,    \cite{RPR10} for
partial sums of weakly
dependent sequences, and  \cite{Fx1,FGL17} for martingales.
\begin{prop}\label{findsa}
Assume condition (\ref{c2}).
Let $p\in (0, 1)$.
Assume that there exist some positive constants  $K_k$  such that, for any $k \in [1, n],$
\begin{equation}\label{laplace2}
 {\mathbb E} \Big[ \big( G_{\varepsilon_k}(\varepsilon_k)\big)^2\exp  \Big(  \big( G_{\varepsilon_k}(\varepsilon_k)\big)^p\Big)\Big] \leq K_k
 \, .
\end{equation}
Set $$K= \sum_{k=1}^{n} \Big(\frac { a_{n-k}(n-k)}{a_{n-1}(n-1)}\Big)^2 K_k.$$
If $K\geq 1,$ then, for any $x \geq0$,
\begin{eqnarray}\label{senib2}
  \mathbb{P}\left(\pm   S_n \geq x   \right)
 \!\!\! &\leq&\!\!\!      2 \exp \Bigg(\frac{- x^2 }{2(K (a_{n-1}(n-1))^2+ x^{2-p}(a_{n-1}(n-1))^p)} \Bigg) .
\end{eqnarray}
\end{prop}

\begin{rem} Let us give some comments on Proposition \ref{findsa}.
\begin{enumerate}
\item If $\sup_k {\mathbb E} \Big[ \big( G_{\varepsilon_k}(\varepsilon_k)\big)^2\exp  \Big(  \big( G_{\varepsilon_k}(\varepsilon_k)\big)^p\Big)\Big] < \infty,$ there exist two positive constants $C$ and $c$ such
that, for any $x >0$,
\begin{eqnarray}\label{jnks}
  \mathbb{P}\left(\pm   S_n \geq n x   \right)
   \leq C   \exp\big( -c\, x^p  n^p   \big),
\end{eqnarray}
where   $C $ and $c$ do not depend on $n$.  For partial sums of independent random variables,
 the rate \eqref{jnks} has been established under semi-exponential moment conditions, see  \cite{LS00}.

\item Notice that $K$ is usually in order of $n.$ Thus the condition $K\geq 1$ is always satisfied for  large enough $n$.
Therefore, inequality  (\ref{senib2}) always holds for large enough $n$.

\item It is worth noting that   for moderate $0\leq  x =o( K^{1/(2-p)})$, the
 bound  (\ref{senib2}) is   sub-Gaussian  and is of the order
\begin{eqnarray}
  \exp\Bigg( -  \frac{x^2}{2K (a_{n-1}(n-1))^2}   \Bigg).
\end{eqnarray}
For  all $x\geq K^{1/(2-p)},$ bound (\ref{senib2}) is a semi-exponential bound and is of the order
\begin{eqnarray}\label{constant}
 \exp\Bigg( - \frac12 \bigg( \frac{x}{a_{n-1}(n-1)} \bigg)^{p}   \Bigg).
\end{eqnarray}
\end{enumerate}
\end{rem}


Without condition (\ref{c2}) but the variables $H_{k,\varepsilon_k}(X_{k-1},\ldots,X_1, \varepsilon_k)$ have semi-exponential moments,   we have the following   semi-exponential deviation inequality.
\begin{prop}\label{propExp}
Let $\alpha \in (0,1)$.
Assume that there exists a constant $C_1$ such that, for any  $k\in [1, n]$,
\begin{equation} \label{Bernsteinmoment}
  \mathbb{E}\,\Big[\exp\Big ( \Big(H_{k,\varepsilon_k}(X_{k-1},\ldots,X_1, \varepsilon_k)\Big)^{\frac{2\alpha}{1-\alpha}} \Big  )  \Big] \leq C_1  .
\end{equation}
Then, for any $x> 0$,
\begin{eqnarray}\label{mainineq01}
  {\mathbb P}\left(\pm  S_n \geq  n x\right)
   \leq  C(\alpha, x)  \exp\left (-\left(\frac{ x  }{ 8 a_{n-1}(n-1)}\right)^{2\alpha} n^\alpha\ \right ) \, ,
\end{eqnarray}
where
\[
C(\alpha,  x)=  2+ 35  C_1     \left( \frac{  a_{n-1}^{2\alpha}(n-1)  }
{ x^{2\alpha} 4^{2-3\alpha}}  + \frac{4 a_{n-1}^2(n-1) }{  x^2 }
\left( \frac{3(1-\alpha)}{2\alpha}\right)^{\frac{ 1-\alpha}{\alpha}} \right)
\]
depends on $n$ only through the term $a_{n-1}(n-1)$.
\end{prop}

\begin{rem}
Clearly, by \eqref{mainineq01}, it holds
\begin{eqnarray}\label{LU00}
\mathbb{P}\left( \pm S_n \geq  n  \right)   = O \Big( \exp\left( - C  n^{\alpha}     \right) \Big)
\end{eqnarray}
for some positive constant $C$.
This convergence rate coincides with the optimal convergence rate for large deviation  of martingales,
as proved in  Theorem 2.1 of  \cite{Fx1} (see also  \cite{LV01} for $\alpha=1/3$).
\end{rem}

\subsection{Fuk-Nagaev type inequalities}\label{Fuksec}

The following proposition follows  from Corollary 2.3 of  \cite{FGL12}.
\begin{prop}\label{po71}
Assume condition (\ref{c2}).
Assume that  there exist some positive constants  $V_k$  such that for any $k \in [1, n],$
$$
{\mathbb E} \big[  \big(   G_{\varepsilon_k}(\varepsilon_k)\big)^2\big] \leq
  V_k  \, .
 $$
Let
\begin{equation}\label{hoeffineq}
V=  \sum_{k=1}^n \Big(a_{n-k}(n-k)\Big)^2 V_k.
\end{equation}
Then, for any $x,y > 0$,
\begin{equation}\label{Hineq}
{\mathbb P}\big(\pm  S_n\geq x \big)  \leq    H_n \left(\frac{x }{y a_{n-1}(n-1)} , \frac{\sqrt{V} }{y a_{n-1}(n-1) }\right)
   + {\mathbb P}\left(  \max_{1 \leq k \leq n} G_{\varepsilon_k}(\varepsilon_k)  > y \right)   ,
\end{equation}
where
\begin{eqnarray}\label{fkmns}
H_{n}(x,v)=\bigg\{\left( \frac{v^2}{
x+v^2}\right)^{ x+v^2 }\left( \frac{n}{n-x}\right)^{ n-x }
\bigg\}^{\frac{n}{n+v^2} }\mathbf{1}_{\{x \leq n\}}
\end{eqnarray}
with the convention that $(+\infty)^0=1$   \emph{(}which  applies when $x=n$\emph{)}.
\end{prop}

\medskip

In particular, when the random variables $(G_{\varepsilon_k}(\varepsilon_k))_{ k\in[1, n]}$ are bounded from above, then Proposition \ref{po71} implies the following Hoeffding bound.
\begin{prop}\label{fhoedding}
Assume condition (\ref{c2}).
Assume that there exist some positive constants $M$  and $V_k$ such that for any $k \in [1, n],$
\begin{equation*}\label{conditionHoeffding}
  \quad G_{\varepsilon_k}(\varepsilon_k) \leq M , \quad \ \
 {\mathbb E}\big[\big(G_{\varepsilon_k}(\varepsilon_k)\big)^2 \big] \leq V_k  .
\end{equation*}
  Then, for any $x > 0$,
\begin{eqnarray}\label{H}
{\mathbb P}\big(\pm S_n>x\big)  \leq     H_n \left(\frac{x }{M a_{n-1}(n-1)} , \frac{\sqrt{V} }{M  a_{n-1}(n-1) }\right) ,
\end{eqnarray}
where $H_{n}(x,v)$ and $V$ are defined by (\ref{fkmns}) and (\ref{hoeffineq}), respectively.
\end{prop}

\begin{rem}
According to Remark 2.1 of \cite{FGL12},  Hoeffding's bound is less than the bounds of Bennett  and Bernstein, that is for any $x\geq0$ and any $v>0$, it holds
\begin{eqnarray}
H_{n}(x, v) \leq  B(x,v):=\left(\frac{v^2}{x+v^2} \right)^{x+v^2}e^x
 \leq  B_1(x,v):= \exp\left(-\frac{x^2}{2(v^2+\frac{1}{3}x )}\right) . \label{Bernstein}
\end{eqnarray}
Then, inequality (\ref{H}) also implies  Bennett's and Bernstein's bounds
\begin{eqnarray*}
\mathbb{P}\big(\pm S_n>x\big)   \leq   B\left(\frac{x }{M a_{n-1}(n-1)} , \frac{\sqrt{V} }{M  a_{n-1}(n-1) }\right)   \leq    B_1\left(\frac{x }{M a_{n-1}(n-1)} , \frac{\sqrt{V} }{M  a_{n-1}(n-1) }\right).
\end{eqnarray*}
\end{rem}

\medskip
For sums of independent random variables with finite moments,
 Fuk \cite{F73} and Nagaev \cite{N79} have established a type of bound combining the Sub-Gaussian bound and polynomial bound.
If the martingale differences $d_i$ have $p$-th moments ($p\geq 2$), then we have the following Fuk-Nagaev type inequality (cf.\ Corollary $3'$ of \cite{F73} and Corollary 1.8 of  \cite{N79}).
\begin{prop}\label{RCfG}
 Assume condition (\ref{c2}).
 Let $p\geq 2$. Assume that there exists some positive constants   $V_k $ and $A_k(p)$ such that
\begin{eqnarray}
{\mathbb E}\big[\big(G_{\varepsilon_k}(\varepsilon_k)\big)^2 \big] \leq V_k\,
 \, \nonumber
\quad   \text{and}  \quad \ \ \,
 {\mathbb E}\big[\big(G_{\varepsilon_k}(\varepsilon_k)\big)^p \big] \leq A_k(p)
  \,  .\label{sdvcx}
\end{eqnarray}
Let $V$ be defined by \eref{hoeffineq}, and let  $$A(p)=  \sum_{k=1}^n\big(a_{n-k}(n-k) \big)^p A_k(p). $$  Then, for any $x  > 0$,
\begin{eqnarray}\label{fuki}
{\mathbb P}( |S_n| \geq x)   \leq  2\Big(1+ \frac2 p \Big)^p \frac{A(p)}{x^p}   + 2 \exp\left( - \frac{2}{(p+2)^2 e^p}  \frac{x^2}{V}  \right).
\end{eqnarray}
\end{prop}
\begin{rem}
Assume that $\sup_{k\geq1}A_k(p) \leq C_p$ for some positive constant $C_p.$
Then $V$ and  $A(p)$ are both of order $n$.  The virtue of bound  (\ref{fuki}) is that it combines the sub-Gaussian bound and polynomial bound
 together. It is easy to see that
the terms
$$\exp\left( - \frac{2}{(p+2)^2 e^p}  \frac{(nx)^2}{V}  \right)\ \ \ \
\textrm{and}\ \ \ \   2\Big(1+ \frac2 p \Big)^p \frac{A(p)}{(xn)^p}$$ are respectively  decreasing at an exponential order  and
at an order $n^{1-p}$.
 Thus, for any $x>0$ and all $n$,
$$
  {\mathbb P}(|S_n|\geq nx)\leq \frac{C_x}{ n^{p-1}}
$$
for some positive $C_x$ not depending on $n$. The last inequality is optimal under the stated condition, even if $S_n$ is a sum of i.i.d.\ random variables.
\end{rem}

\subsection{McDiarmid inequality}\label{McDsec}
Where the increments $d_k$ are bounded, we shall use an improved version of the well
known inequality by McDiarmid, which has been established by \cite{R13}.
Following the notations in   \cite{R13},  we denote
$$\ell(t)= (t - \ln t -1) +t (e^t-1)^{-1} + \ln(1-e^{-t}), \ \ \ \ \
\ \ t > 0, $$ and let
$$\ell^*(x)=\sup_{t>0}\big(xt- \ell (t) \big), \ \ \ \ \  \ \ x > 0, $$ be the Young transform of $\ell(t)$.
 It is known that for any $x \in [0,1)$, the following  inequalities hold
\begin{equation}\label{is50}
\ell^*(x) \geq (x^2-2x) \ln(1-x) \geq 2x^2 + x^4/6\, ,
\end{equation}
see \cite{R13}.
Let   $(\varepsilon'_i)_{i \geq 1}$ be an independent copy of
$(\varepsilon_i)_{i \geq 1}$.
\begin{prop} \label{classic}
Assume that there exist some positive constants $M_k$ such that
\begin{equation}\label{in850}
\big \|d\big(F_k ( (X_{k-i})_{i\geq 1}   ;   \varepsilon_k), F_k ((X_{k-i})_{i\geq 1}  ; \varepsilon'_k)\big)
\big \|_\infty \leq M_k, \ \   \text{ $k \in [1,   n] $.}
\end{equation}
Let
$$
M^2(n,\rho)= \sum_{k=1}^n \big(a_{n-k}(n-k) M_k \big)^2 \quad \text{and} \quad D (n, \rho)=\sum_{k=1}^n a_{n-k}(n-k) M_k \, .
$$
Then, for any $t \geq 0$,
\begin{eqnarray}\label{rio1}
\mathbb{E}[e^{\pm   tS_n } ]  \leq   \exp \left (  \frac{D^2(n, \rho)}{M^2(n,\rho)}\  \ell \Big( \frac {M^2(n,\rho)\, t} {D(n, \rho)} \Big)\right) \,
\end{eqnarray}
and,  for any   $x \in [0, D(n, \rho)]$,
\begin{eqnarray}\label{rio2}
{\mathbb P}\big(\pm  S_n>x\big) \leq \exp \left ( -\frac{D^2(n, \rho)}{M^2(n,\rho)} \ \ell^*\Big( \frac {x} {D(n, \rho)} \Big)\right) \, .
\end{eqnarray}
Consequently, for any $x \in [0, D(n, \rho)]$,
\begin{eqnarray}\label{riosbound}
{\mathbb P}\big(\pm  S_n>x\big)  \leq  \left ( \frac{D(n, \rho)-x}{D(n, \rho)}\right)^{ \frac{2D(n, \rho)x -x^2}{M^2(n,\rho)}} .
\end{eqnarray}
\end{prop}

\begin{rem} Since for any $x \in [0,1)$, $(x^2-2x) \ln(1-x) \geq 2\,x^2$, inequality
(\ref{riosbound}) implies the following McDiarmid inequality
\begin{eqnarray*}
{\mathbb P}\big(\pm  S_n>x\big)   \leq    \exp \left ( -\frac{2 x^2}{M^2(n,\rho)}
 \right) \, .
\end{eqnarray*}
\end{rem}

\begin{rem}
Taking $\Delta (n, \rho)=a_{n-1}(n-1)\max_{1\leq k \leq n} M_k$, we obtain the upper bound: for any $x \in [0, n \Delta(n, \rho)]$,
\begin{eqnarray*}\label{rio2bis}
{\mathbb P}\big(\pm  S_n>x\big)
 \leq \exp \left ( -n \ell^*\Big( \frac {x} {n \Delta(n, \rho)} \Big)\right) \leq \exp \left ( -\frac{2 x^2}{n\Delta^2(n,\rho)}
 \right)
  \, .
\end{eqnarray*}
\end{rem}

\subsection{von Bahr-Esseen bound}\label{VBEB}
In the first proposition of this section, we assume that the dominating random variables
$G_{\varepsilon_k}(\varepsilon_k)$ have only a moment of order $p \in [1,2]$.
For similar inequalities in the case where the
$X_i$'s are independent, we refer to   \cite{P10}.
\begin{prop}\label{VBEI}
Let $p \in [1, 2]$. Assume that for any $k \in [1, n],$
\begin{equation}\label{fnf58}
{\mathbb E} \Big[  \Big(  H_{k, \varepsilon_k}(X_{k-1},\ldots,X_1, \varepsilon_k)\Big)^p\Big] \leq A_k(p) \, .
\end{equation}
Then
\begin{equation}\label{vBE}
  \| S_n \|_p  \leq    \big( A(n, p) \big)^{1/p},
\end{equation}
 where
\begin{equation}\label{Ap}
A(n, p)= A_1(p)\big(a_{n-1}(n-1)\big)^p    +
2^{2-p} \sum_{k=2}^n \big(a_{n-k}(n-k)\big)^p A_k(p) \, .
\end{equation}
By (\ref{vBE}), it follows that, for any $x  > 0$,
\begin{eqnarray} \label{fdgsd}
{\mathbb P}( |S_n| \geq nx)   \leq  \frac{A(n, p)}{(nx)^p }  .
\end{eqnarray}
\end{prop}
\begin{rem}
The constant $2^{2-p}$ in \eref{Ap} can be improved to the more precise constant $\tilde C_p$ described in Proposition 1.8 of   \cite{P10}.
\end{rem}
\begin{rem} If $\sup_k A_k(p) < \infty,$ then (\ref{fdgsd}) implies that, for any $x>0$ and all $n$,
$$
  {\mathbb P}(|S_n|\geq nx)\leq \frac{C_x}{ n^{p-1}}
$$
for some positive $C_x$ not depending on $n$.
\end{rem}

For any real-valued random variable $Z$ and any $p\geq 1$, define the weak
moment of order $p$  by
\begin{equation}\label{weakp}
\|Z\|_{w,p}^p=\sup_{x>0} x^p{\mathbb P}(|Z|>x)\, .
\end{equation}
We now consider the case where the variable  $  H_{k, \varepsilon_k}(X_{k-1},\ldots,X_1, \varepsilon_k)$
has only a weak moment of order $p \in (1,2)$.

\begin{prop}\label{weakVBEI}
Let $p \in (1, 2)$. Assume that for any $k \in [1, n],$
\begin{equation}\label{weakVB}
\Big \|   H_{k, \varepsilon_k}(X_{k-1},\ldots,X_1, \varepsilon_k)\Big \|_{w,p}^p \leq A_k(p)  .
\end{equation}
Then, for any $x>0,$
\begin{equation}\label{weakVBE}
 {\mathbb P}(|S_n| \geq x) \leq  \frac{C_p B(n,p)}{x^p}\, ,
\end{equation}
where
 $$
 C_p= \frac{4p}{(p-1)} +\frac{8 }{(2-p)}\,
\ \ \ \
\textrm{and} \ \ \ \
B(n, p)=
 \sum_{k=1}^n \big( a_{n-k}(n-k)\big)^p A_k(p) \,.
$$
\end{prop}

\begin{rem} If $\sup_{k\geq1}A_k(p) < \infty,$ then, from (\ref{weakVBE}), we deduce that, for any $x>0$ and all $n$,
$$
  {\mathbb P}(|S_n| \geq nx)\leq \frac{C_x}{ n^{p-1}} \, ,
$$
for some positive $C_x$ not depending on $n$. The last inequality is optimal under the stated condition, even if $S_n$ is a sum of i.i.d.\ random variables.
\end{rem}

\section{Moment inequalities for the functional $S_n$}\label{section5}
\setcounter{equation}{0}

In this section, we also assume that $(X_{i})_{i\leq 0}$ are deterministic.
We present some moment inequalities for the functional $S_n$.

\subsection{Marcinkiewicz-Zygmund bound}\label{MZB}
Now assume that the  random variables
$( H_{k,\varepsilon_k}(X_{k-1},\ldots,X_1, \varepsilon_k))_{k \in [1, n]}$ have  moments of order $p \geq 2$.
\begin{prop}\label{MZIP}
Let $p \geq 2$. Assume that for any $k \in [1, n],$
\begin{equation}\label{shnbs}
{\mathbb E} \Big[  \Big(   H_{k,\varepsilon_k}(X_{k-1},\ldots,X_1, \varepsilon_k)\Big)^p\Big] \leq A_k(p)\, .
\end{equation}
Then
\begin{equation}\label{MZI}
  \| S_n \|_p \leq   \sqrt{ A(n,p)}\, ,
\end{equation}
 where
 $$
A(n, p)= \big(a_{n-1}(n-1) \big)^2 \big(A_1(p) \big)^{2/p} +
(p-1) \sum_{k=2}^n \big(a_{n-k}(n-k) \big)^2\big(A_k(p) \big)^{2/p} \, .
$$
\end{prop}
\begin{rem}
Assume that  $F_k$ satisfies   \eref{contract} and (\ref{c2}). Then, it follows from the proof of Proposition \ref{MZIP} that
the inequality
\eref{MZI} remains true if  condition   \eref{shnbs} is replaced by, for any $k \in [1, n],$
\begin{equation*}
  {\mathbb E} \big [ \big ( G_{\varepsilon_k}(  \varepsilon_k) \big)^p
 \big ] \leq A_k(p) \, .
\end{equation*}
\end{rem}
\begin{rem}
Assume that $\sup_{k\geq1}A_k(p) \leq C_p$ for some positive constant $C_p.$ Then, from (\ref{MZI}), we deduce that
\begin{equation}
  \| S_n \|_p =O(\sqrt{n}\,  ), \ \ n\rightarrow \infty.
\end{equation}
Moreover, by Markov's inequality, we have,  for  any $x>0,$
\begin{equation}
 {\mathbb P}(|S_{ n}|\geq x ) \leq  \frac{  C_{1,p}}{x^p}\, \sqrt{ A(n, p)}
\end{equation}
and
$$
  {\mathbb P}(|S_n| \geq nx)\leq \frac{C_{x, p}}{ n^{p/2}} \, ,
$$
for some positive $C_{x, p}$ not depending on $n$.
\end{rem}

\noindent\emph{Proof.} By Theorem 2.1 of  \cite{R09},  we get, for $p \geq 2,$
 $$
\| S_n \|_p^2  \leq \|d_1\|_p^2 +(p-1) \sum_{k=2}^n  \| d_k \|_p^2 \, .
$$
By point 2 of Proposition \ref{McD}  and condition (\ref{shnbs}), we can deduce that
\begin{eqnarray*}
\| S_n \|_p^2   &\leq&   \big( a_{n-1}(n-1) \big)^2  \left( {\mathbb E} \big[  \big( H_{1,\varepsilon_1}( \varepsilon_1)\big)^p\big]  \right)^{2/p} \\
&&  + (p-1)
    \sum_{k=2}^n \big(  a_{n-k}(n-k) \big)^2 \left({\mathbb E} \Big[  \Big(  H_{k,\varepsilon_k}(X_{k-1},...,X_1, \varepsilon_k)\Big)^p\Big]\right)^{2/p} \\
    &\leq&  A(n, p),
\end{eqnarray*}
which gives the desired inequality. \hfill\qed

\subsection{Rosenthal bounds}\label{BRB}
When the random variables $(G_{\varepsilon_k}(\varepsilon_k)  )_{1\leq k \leq n}$  have   moments of order $p \geq 2$,
we have the following proposition. For similar inequalities for sum of independent random variables, we refer to \cite{P13}.
\begin{prop}\label{Ros}
Assume condition (\ref{c2}) holds.
Assume that there exist some  constants  $ V_k \geq 0$   such that
\begin{equation}\label{Strongros}
 {\mathbb E}\big[\big( G_{\varepsilon_k}(\varepsilon_k)\big)^2 \big] \leq V_k
\, .
\end{equation}
Let
\begin{equation}\label{variance}
V=  \sum_{k=1}^n \big(a_{n-k}(n-k) \big)^2 V_k\, .
\end{equation}
For any  $p\geq 2$, there exist two positive constants $C_{1,p}$ and $C_{2,p}$ such that
\begin{equation}\label{rosI}
  \| S_n \|_p \leq   C_{1,p}  \sqrt{ V} + C_{2,p}  \Big \|  \max_{1 \leq k \leq n} a_{n-k}(n-k)G_{\varepsilon_k}(\varepsilon_k) \Big \|_p  .
\end{equation}
\end{prop}

\noindent\emph{Proof.}
Applying Proposition \ref{McD}, we have  $|d_k|\leq a_{n-k}(n-k) G_{\varepsilon_k}(\varepsilon_k)$ for any $k \in [1,  n],$
and consequently
$$
 \quad {\mathbb E}[d_k^2|{\mathcal F}_{k-1}] \leq  \big(a_{n-k}(n-k)\big)^2 V_2\, \ \ \ \text{ for any $k\in [1,  n].$}
$$
Then the proposition  follows  directly from Theorem 4.1 of  \cite{P94}. \hfill \qed

\begin{rem}
According to the proof of Theorem 4.1 of  \cite{P94}, one can take $C_{1,p} =60c$ and $C_{2,p} = 120\sqrt{c} \,e^{p/c}$ for any $c \in [1,p]$.
\end{rem}

\medskip

Next, we consider the case where the random variables $(G_{\varepsilon_k}(\varepsilon_k)  )_{1\leq k \leq n}$  have  a weak moment of order $p \geq2$.
Recall that the weak moment $\| Z \|_{w,p}^p$ has been defined by \eref{weakp}.
\begin{prop}\label{RCG}
Assume that the conditions (\ref{c2}) and \eref{Strongros} hold, and let $V$ be defined by \eref{variance}.
Then, for any  $p\geq 2$, there exist two positive constants $C_{1,p}$ and $C_{2,p}$ such that  for any $x > 0,$
\begin{equation}\label{weakrosI}
  {\mathbb P}(|S_n|\geq x) \leq
  \frac{1}{x^p} \left ( C_{1,p} {V}^{p/2} + C_{2,p}  \Big \|   \max_{1 \leq k \leq n} a_{n-k}(n-k)G_{\varepsilon_k}(\varepsilon_k) \Big \|^p_{w,p} \right ).
\end{equation}
\end{prop}

\noindent\emph{Proof.} The argument is similar to the proof  of Proposition \ref{Ros}, by applying Theorem 6.3 in \cite{CG12}. \hfill \qed

\section{Applications}
\setcounter{equation}{0}

\subsection{Application to  stochastic gradient Langevin dynamic}\label{SGLD}
Let $\psi: \mathbb{R}^d \times \mathbb{R}^r \longrightarrow \mathbb{R},   (\omega,\zeta) \mapsto \psi(\omega,\zeta) $, be a non-convex stochastic loss function. Consider the optimization problem
\begin{eqnarray*}
	\omega^*=\text{argmin}_{\omega\in\R^d}P(\omega),\quad P(\omega)=\E_{\zeta\sim\nu}\psi(\omega,\zeta),
\end{eqnarray*}
where $\zeta$ is a random variable with probability distribution $\nu$. To find the optimizer $\omega^*$, Welling and Teh \cite{welling2011bayesian} proposed the stochastic gradient Langevin dynamic (SGLD) algorithm: for $k\geq 1,$
\begin{eqnarray}\label{e:sgld}
	\omega_k=\omega_{k-1}-\eta\nabla\psi(\omega_{k-1},\zeta_k)+\sqrt{\eta\delta}\xi_{k},
\end{eqnarray}	
where $\eta>0$ is the step size, $\delta>0$ is the inverse temperature parameter, $(\xi_k)_{k\ge1}$ is a  sequence of  i.i.d.\   random vectors with $\cov[\xi_{k} ]=   I_d,$ where $I_d$ is an identity $d$-dimensional matrix, and $(\zeta_k)_{k\ge1}$ are i.i.d.\ samples from $\nu$.
Rewriting \eqref{e:sgld}, we have
\begin{eqnarray}\label{e:sgld1}
	\omega_k&=&\omega_{k-1}-\eta\nabla P(\omega_{k-1})+\eta\nabla P(\omega_{k-1})-\eta\nabla\psi(\omega_{k-1},\zeta_k)+\sqrt{\eta\delta}\xi_{k}\nonumber\\
	&:=&\omega_{k-1}-\eta\nabla P(\omega_{k-1})+\sqrt\eta V_{\eta,\delta}(\omega_{k-1},\zeta_k,\xi_{k}),
\end{eqnarray}
where
\begin{eqnarray*}
	V_{\eta,\delta}(\omega_{k-1},\zeta_k,\xi_{k})&=& \sqrt\eta\nabla P(\omega_{k-1})-\sqrt\eta\nabla\psi(\omega_{k-1},\zeta_k)+\sqrt{\delta}\xi_{k}.
\end{eqnarray*}
As  $\E\psi(\cdot,\zeta)=P(\cdot)$, by some simple calculations, we have
\begin{eqnarray*}
	\E[V_{\eta,\delta}(\omega_{k-1},\zeta_k,\xi_{k})\ |\   \omega_{k-1}  ] = 0
\end{eqnarray*}
and
\begin{eqnarray*}
\cov[V_{\eta,\delta}(\omega_{k-1},\zeta_k,\xi_{k})\,|\,\omega_{k-1}] &=  &   \E[V_{\eta,\delta}(\omega_{k-1},\zeta_k,\xi_{k})V_{\eta,\delta}(\omega_{k-1},\zeta_k,\xi_{k})^\top \,|\, \omega_{k-1}]\\
 & = & \eta\Sigma(\omega_{k-1})+\delta I_d,
\end{eqnarray*}
where
\begin{eqnarray*}
	\Sigma(x)=\E[\nabla\psi(x,\zeta)\nabla\psi(x,\zeta)^\top]-\nabla P(x)\nabla P(x)^\top.
\end{eqnarray*}
 For the cost function $\psi$ and random variable $\zeta$, we introduce the following conditions. Assume that there exist  constants $L, M, K >0$  such that for any $x,y, z\in \mathbb{R}^d$,
	\begin{eqnarray}\label{e:lip}
		\|\nabla\psi(x,z)-\nabla\psi(y,z) \|_2\le L \|x-y \|_2,
	\end{eqnarray}
\begin{eqnarray}\label{e:lip12}
		\|\nabla\psi(x,y)-\nabla\psi(x,z)\|_2\le M \|y-z\|_2,
	\end{eqnarray}
	\begin{eqnarray}\label{e:dissi}
		\langle x-y, \nabla \psi(x,z)-\nabla \psi(y,z) \rangle \geq  K \|x-y\|_2^2.
	\end{eqnarray}
Then, by \eref{e:sgld}, the condition (\ref{contract}) is satisfied with $$F_n(   (x_{n-i})_{i\geq 1} ;   \varepsilon_n)=   x_{n-1}-\eta\nabla\psi(x_{n-1},\zeta_n)+\sqrt{\eta\delta}\xi_{n} ,\, \ \ \ \  \varepsilon_n^T=(\zeta_n, \xi_{n}),  $$
$a_1=(1- 2\eta K + \eta^2L^2 )^{1/2}$ and $a_i=0, i\geq 2$, provided that $0<\eta  <\{1/2K, 2K/L^2  \}$.
Indeed, by the conditions (\ref{e:lip})  and (\ref{e:dissi}), it is easy to see that
\begin{eqnarray*}
{\mathbb E}\|F(x; \varepsilon_n  )-F(x';  \varepsilon_n )\|_2^2 &=& {\mathbb E}\|x-x'-\eta\nabla\psi(x,\zeta_n) + \eta\nabla\psi(x',\zeta_n) \|_2 \\
&=&\|x-x'\|_2+ \eta^2{\mathbb E}\|  \nabla\psi(x,\zeta_n) - \nabla\psi(x',\zeta_n) \|_2^2 \\
&& - 2\eta (x-x')^T( \nabla\psi(x,\zeta_n) - \nabla\psi(x',\zeta_n))\\
&\leq&(1- 2\eta K + \eta^2L^2 )\|x-x'\|_2^2,
\end{eqnarray*}
which implies that
\begin{eqnarray*}
{\mathbb E}\|F(x; \varepsilon_n  )-F(x';  \varepsilon_n )\|_2   \leq
  \sqrt{1- 2\eta K + \eta^2L^2} \|x-x'\|_2 .
\end{eqnarray*}
Thus, the condition (\ref{contract}) is satisfied.
Moreover,   by the conditions  (\ref{e:lip12}) and (\ref{e:dissi}),  it holds
\begin{eqnarray*}
 \|F(x; \varepsilon_n  )-F(x;  \varepsilon_n' )\|_2^2 &=&  \| \eta\nabla\psi(x ,\zeta_n)+\sqrt{\eta\delta}\xi_{n} -
\eta\nabla\psi(x ,\zeta_n')-\sqrt{\eta\delta}\xi_{n}' \|_2^2 \\
&\leq& 2 \eta^2 \|\nabla\psi(x ,\zeta_n)- \nabla\psi(x ,\zeta_n')\|_2^2 + 2\eta\delta\| \xi_{n} - \xi_{n}' \|_2 ^2  \\
&\leq& 2 \eta^2 M^2 \| \zeta_n -  \zeta_n' \|_2^2 + 2\eta\delta\| \xi_{n} - \xi_{n}' \|_2 ^2  \\
&\leq& 2 (\eta^2 M^2 \vee \eta\delta) \|    \varepsilon_n    -    \varepsilon_n' \|_2^2 ,
\end{eqnarray*}
which implies that
\begin{eqnarray*}
 \|F(x; \varepsilon_n  )-F(x;  \varepsilon_n' )\|_2
 \leq  \sqrt{2 (\eta^2 M^2 \vee \eta\delta)} \|    \varepsilon_n    -    \varepsilon_n' \|_2 .
\end{eqnarray*}
Thus, the condition (\ref{c2}) is  also satisfied with
$$
G_{\varepsilon_k}(y)= \sqrt{2 (\eta^2 M^2 \vee \eta\delta)} \int  \|    y    -    x   \|_2   P_{\varepsilon_k}(dx ) \, .
$$
Assume that  \eqref{e:sgld} is exponential ergodic with invariant measures $\pi_\eta$. As $\omega_k$ weak convergences to $\pi_\eta$,   define
$$  \Pi_\eta(\cdot)=  \frac{1}{[\eta^{-2}]} \sum_{k=0}^{[\eta^{-2}]-1  } \delta_{\omega_k}(\cdot) ,$$
where $\delta_{y}(\cdot)$ is the Dirac measure of $y.$ Here, we denote $[\eta^{-2}]$   the largest integer less than $  \eta^{-2} $.

We consider the asymptotic property of $\Pi_\eta(h)$,
where $h: \R^d \rightarrow \R$ is a $1$-Lipschitz test function.
Assume that the initial value $ \omega_0$ of SGLD  algorithm is  deterministic. As an illustration for our results, assume that $\delta$ is a positive constant. Denote
$$    \hat{\varepsilon}_k    =  \int  \|     \varepsilon_k   -    x   \|_2   P_{\varepsilon_k}(dx ) .$$
We have the following qualitative consequences of these inequalities in Sections \ref{deviationiq} and \ref{section5}:
\begin{itemize}
\item If there exists a constant $M$  such that for  any $l\geq 2$,
\begin{equation}
 {\mathbb E} \,  \hat{\varepsilon}_k ^{\,l}  \leq
 \frac {l!}{2}   M^{l-2}   {\mathbb E} \, \hat{\varepsilon}_k^{\,2}  ,
\end{equation}  then
there exist some positive constants $A$ and $B$ such that
\begin{equation}\label{LW}
{\mathbb P}\bigg(\big| \Pi_\eta(h)  - \pi_\eta(h)  \big|>x\bigg) \leq \begin{cases}
  2\exp\Big(- \eta^{-3/2}  A x   \Big) \quad \ \  \, \text{if\ \ $x \geq \eta^{-1/2}, $}\\
  \\
   2\exp\Big(-  \eta^{-1}  B x^2  \Big) \quad  \ \ \  \, \text{if\ \ $x \in [0, \eta^{-1/2}]$.}
  \end{cases}
\end{equation}
This follows from Proposition \ref{pr01}.

\item
Let $p\in (0, 1)$. If there exist some positive constants $K$ and $L$ such that
\begin{equation}\label{laplace2}
 {\mathbb E} \Big[   \hat{\varepsilon}_k ^{\,2}\exp  \Big(   L \, \hat{\varepsilon}_k^{\,p}\Big)\Big] \leq K
 \, ,
\end{equation}
 then there exist some positive constants $A, B, C $ and $D $ such that
$$
{\mathbb P}\bigg(\big|   \Pi_\eta(h)  - \pi_\eta(h)   \big|>x\bigg) \leq \begin{cases}
  C\exp\Big(- \eta^{- 3p/2}   A x^{p}   \Big) \quad \ \ \  \ \text{if\ \ $x \geq   \eta^{  (1-3p/2)/(2-p)},$}\\
  \\
   D\exp\Big(- \eta^{-1}  B x^2  \Big) \quad \ \ \ \ \ \ \  \text{if\ \ $x \in [0,   \eta^{ (1-3p/2)/(2-p)}]$.}
  \end{cases}
$$
This follows from Proposition \ref{findsa}.

\item Let $p\geq 2$. If
\begin{equation}\label{fnf5ds8}
   {\mathbb E}\, \hat{\varepsilon}_k^{\,p}   < \infty     ,
\end{equation}
  then there exists a positive constant $C$
such that
$$
{\mathbb P}\Big(\big|   \Pi_\eta(h)  - \pi_\eta(h)    \big|>x\Big) \leq \frac{C}{
  x^p} \eta^{ p/2 } \, .
$$
This follows from Proposition \ref{RCG}.

\item If \eref{fnf5ds8} holds for some  $p \in (4/3, 2)$,
  then there exists a positive constant $C$
such that
$$
{\mathbb P}\Big(\big|   \Pi_\eta(h)  - \pi_\eta(h)    \big|>x\Big) \leq \frac{C}{
  x^p} \eta^{3p/2-2} \, .
$$
This follows from   Proposition \ref{weakVBEI}.
Notice that in this model, we have
\begin{equation}\label{lapdgl56}
  H_{k,\varepsilon_k}(X_{k-1},\ldots,X_1, \varepsilon_k) \leq \hat{\varepsilon}_k\ \ \ \textrm{and} \ \ \ \big \|  \hat{\varepsilon}_k  \big \|^p_{w,p}  \leq   \| \hat{\varepsilon}_k   \|_{p} .
\end{equation}
\end{itemize}
Using \eqref{lapdgl56}, we have the following  moment bounds for $ \Pi_\eta(h)  - \pi_\eta(h)$:
\begin{itemize}

\item  If \eref{fnf5ds8} holds for some $p \geq2$, then
\begin{equation} \nonumber
\Big\|  \Pi_\eta(h)  - \pi_\eta(h)   \Big\|^p_p \leq C  \, \eta^{p/2}.
\end{equation}
This follows from  \eqref{MZI}.

\item  If \eref{fnf5ds8} holds for some $p \in (4/3, 2)$, then
\begin{equation} \nonumber
\Big\|  \Pi_\eta(h)  - \pi_\eta(h)   \Big\|^p_p \leq C  \, \eta^{3p/2-2}.
\end{equation}
This follows from    \eqref{vBE}.
\end{itemize}

Let us now give some references on the subject. By \eqref{e:sgld1},
it is nature to consider the following SDE to approximate \eqref{e:sgld}, that is
\begin{eqnarray}\label{sddsfsd}
	d X_t=-\nabla P(X_t)d t+Q_{\eta,\delta}(X_t)d B_t,
\end{eqnarray}
where $Q_{\eta,\delta}(x)=\big(\eta\Sigma(x)+\delta I_d\big)^{\frac12}$
is a positive definite matrix and $B_t$ is a $d-$dimensional standard Brownian motion.
When $(\xi_k)_{k\ge1}$ is a sequence of i.i.d.\ standard $d$-dimensional normal random vectors and  the random variable $\nabla\psi(x,\zeta)$ is sub-Gaussian for any $x \in \mathbb{R}^d $, that is, there exist positive constants $C_\zeta$ and $C$ such that
\begin{eqnarray}
\E\exp\{C_\zeta|\nabla\psi(x,\zeta)|^2\} \le C.
\end{eqnarray}
 Dai et al.\ \cite{DFL25} proved that
\begin{eqnarray*}
	W_1(\pi,\pi_\eta)=\sup_{h\in\mathrm{Lip_1}}|\pi_\eta(h)-\pi(h)|=O(\eta^{1/2}),
\end{eqnarray*}
where $\pi$ is the invariant measure of the SDE (\ref{sddsfsd}).

\subsection{Conclusions}
The deviation inequalities certainly have a lot of applications. We refer to  \cite{DF15} for the convergence rates in the Wasserstein distance between the empirical distribution and the invariant distribution.
Applications to empirical risk minimization and stochastic approximation by averaging for linear problem can also be found
in  \cite{FAD22}. It is also easy to see that the deviation inequalities are applicable to  mean fields memory models in Subsection \ref{secexample}.

\paragraph{\bf Acknowledgment.}
Xiequan Fan was partially supported by the National Natural Science Foundation
of China (Grant Nos.\,12371155 and 11971063). Paul Doukhan was also funded by CY-AS
(``Investissements d'Avenir" ANR-16-IDEX-0008), ``EcoDep" PSI-AAP2020-0000000013.
The first author is also thankful for the very warm and kind support of Northeastern University at Qinhuandao.

\newpage

\section*{Appendix}\label{appendix}
\setcounter{equation}{0}

In this appendix, we present the proof  of Proposition \ref{McD} and proofs of the propositions in Section \ref{deviationiq}.

\vspace{0.3cm}

\noindent {\emph{Proof of Proposition \ref{McD}.}} The first point will be proved by recurrence in the backward sense.
The result is obvious for $k=n$, since $g_n=f$. For $k=n-1,$ it holds
\begin{multline}
g_{n-1}(X_1, X_2,\ldots, X_{n-1})={\mathbb E}[g_n(X_1, X_2, \ldots, X_n)|{\mathcal F}_{n-1}]  \\
 \ \ \ \ = \int\!\!\!\!\int g_n(X_{1}, X_2,\ldots, X_{n-1}, F_n( (X_{n-1},\ldots,X_1,x_0,\ldots) ;y)) P_{\varepsilon_n} ( dy)\widetilde{P} ( dx_0,\ldots)\, .
\end{multline}
Set $x_i=x_i'$ for $i\leq 0.$
It is easy to see that
\begin{multline}\label{triv1}
|g_{n-1}(x_1, x_2, \ldots, x_{n-1})-g_{n-1}(x'_1, x'_2, \ldots, x'_{n-1})|\\ \leq \int\!\!\!\!\int \big|g_n(x_1,\ldots,x_{n-1},  F_n((x_{n-i})_{i\geq 1};y))-g_n(x'_{1},\ldots, x'_{n-1}, F_n((x_{n-i}')_{i\geq 1};y)) \big|
P_{\varepsilon_n}(dy)\,\widetilde{P} ( dx_0,\ldots)  \\
 \leq   d(x_1,x'_1)+\cdots +    d(x_{n-1}, x'_{n-1})   +   \int \!\!\!\!  \int  d(F_n( (x_{n-i})_{i\geq 1};y), F_n((x'_{n-i})_{i\geq 1};y))P_{\varepsilon_n}(dy)\,\widetilde{P} ( dx_0,\ldots) \\
 \leq   d(x_1,x'_1)+\cdots +   d(x_{n-1}, x'_{n-1})   +    \sum_{i=1}^{n-1} a_{i}  \,  d(x_{n-i}, x_{n-i}')   \\
 \leq  \sum_{i=1}^{n-1}(1+ a_{i}  )  d(x_{n-i}, x'_{n-i})  .\ \ \quad \quad  \quad  \quad \quad \quad  \quad  \quad \quad \quad  \quad  \quad
\end{multline}
Set
\begin{eqnarray}\label{fsjdsfm}
  a_1(i)=1+a_i,  \ \  a_{k+1}(i)=a_{k}(i)+ a_{k}(k) a_{i-k}, \ \ \  k \in [1, \, n-1]\  \textrm{and} \ i \in [k+1, n-1].
\end{eqnarray}
Assume
\begin{eqnarray}\label{fsjkm}
|g_{n-k}(x_1, x_2,  \ldots, x_{n-k})-g_{n-k}(x'_1, x'_2, \ldots, x'_{n-k})| \leq  \sum_{i=k}^{n-1} a_{k}(i) \, d(x_{n-i}, x'_{n-i})
\end{eqnarray}
holds for $k=n-k.$
Then for $k=n-(k+1),$ we have
\begin{multline}
|g_{n-k-1}(x_1, x_2, \ldots, x_{n-k-1})-g_{n-k-1}(x'_1, x'_2, \ldots, x'_{n-k-1})|\\ \leq  \int \!\!\!\!\int \big|g_{n-k}(x_1,\ldots,x_{n-k-1},  F_{n-k }((x_{n-i})_{i\geq k+1};y))  \quad \quad \quad \quad \quad \quad \quad \quad \quad \quad \quad \quad \quad \quad \quad \quad  \\ \quad \quad \quad \quad \quad \quad   -g_{n-k}(x'_{1},\ldots, x'_{n-k-1}, F_{n-k }( (x'_{n-i})_{i\geq k+1};y)) \big|
P_{\varepsilon_{n-k}}(dy)\,  \widetilde{P} ( dx_0,\ldots\ldots) \\
 \leq  \sum_{i=k+1}^{n-1} a_{k }(i)  d(x_{n-i}, x'_{n-i})  \quad \quad  \quad \quad \quad \quad \quad \quad \quad \quad \quad \quad \quad \quad \quad \quad \quad \quad \quad \quad \quad \quad \quad \quad \\
  \quad \quad    +     a_{k}(k) \int\!\!\!\!\int d(F_{n-k }((x_{n-i})_{i\geq k+1};y), F_{n-k }((x'_{n-i})_{i\geq k+1};y))P_{\varepsilon_{n-k}}(dy)\,  \widetilde{P} ( dx_0,\ldots) \\
 \leq  \sum_{i=k+1}^{n-1} a_{k }(i)  d(x_{n-i}, x'_{n-i})  +     a_{k}(k) \sum_{i=k+1}^{n-1}  a_{i-k}\, d(x_{n-i}, x'_{n-i})   \quad \quad \quad \\
 \leq   \sum_{i=k+1}^{n-1}  a_{k+1}(i)  d(x_{n-i}, x'_{n-i})  , \ \ \   \ \ \ \ \ \   \ \ \
\end{multline}
which justifies that inequality (\ref{fsjkm}) holds for each $k \in [1, n]$.
Clearly, inequality (\ref{fsjkm}) is equivalent to the first desired inequality of Proposition  \ref{McD}.
Using   equality (\ref{fsjdsfm}) and the definition of $ a_{k}(i)$, we deduce that
\begin{eqnarray}
  a_{k+1}(k+1)&=&a_{k }(k+1)+ a_{k }(k ) a_{1} \ =\  a_{k-1 }(k+1)+  a_{k -1}(k-1 ) a_{2}+ a_{k }(k ) a_{1} \nonumber \\
  &=& a_1(k+1)+ \sum_{l=1}^{k }  a_{l}(l) a_{k+1-l}.   \nonumber
\end{eqnarray}
Since $a_1(i)= 1+a_i, $  we obtain
\begin{eqnarray}
  a_k(k)
   = 1+a_k+    \sum_{l=1}^{k-1}  a_{l}(l) a_{k -l}.
\end{eqnarray}
This completes the proof of the  point 1.

Let us prove the point 2.
In the same way, for $k=1,$
\begin{align*}
 |d_1| & = \big|g_1(X_1)-{\mathbb E}[g_1(X_1)]\big|\
  \leq  \int   \big |g_1( X_1 )     -g_1( y) \big| P_{X_1}(dy) \\
 &\leq a_{ n-1}(n-1)  \int   d( X_1,
 y)P_{X_1}(dy)   \ = \ a_{ n-1}(n-1)
  H_{1, X_1}(X_1) \,,
\end{align*}
and for any $k \in [2,n]$,
\begin{align*}
 |d_k| & = \big|g_k(X_1, \cdots, X_k)-{\mathbb E}[g_k(X_1, \cdots, X_k)|{\mathcal F}_{k-1}]\big|\\
 &\leq  \int\!\!\!\! \int \big |g_k(X_1, \cdots, X_{k-1}, F_k(\mathbf{X}_{k-1},x_0,\ldots; \varepsilon_k)) \\
 & \ \ \ \ \ \ \ \ \ \ \ \  \ \   -g_k(X_1, \cdots,X_{k-1}, F_k(\mathbf{X}_{k-1},x_0,\ldots; y))\big| P_{\varepsilon_k}(dy)\widetilde{P} ( dx_0\ldots)\\
 &\leq a_{ n-k}(n-k)  \int\!\!\!\!\int  d(F_k(\mathbf{X}_{k-1},x_0,\ldots; \varepsilon_k),
 F_k(\mathbf{X}_{k-1},x_0,\ldots; y))P_{\varepsilon_k}(dy)\widetilde{P} ( dx_0\ldots) \\
 &  = a_{ n-k}(n-k)
  H_{k,\varepsilon_k}(\mathbf{X}_{k-1}, \varepsilon_k) \,,
\end{align*}
which completes the proof of point 2.

The  point 3 follows easily from  point 2, by using (\ref{c2}).

For point 4, we only need to prove it for $k=1$ and the remaining follows by point 3.
As $(X_{i})_{i\leq 0}$ are deterministic, we can deduce  that
\begin{align*}
 |d_1| & = \big|g_1(X_1)-{\mathbb E}[g_1(X_1)]\big|
  \ \leq \ \int  \big |g_1( X_1 )     -g_1(  F_1( X_0,\ldots; y)) \big| P_{\varepsilon_1}(dy) \\
 &\leq a_{ n-1}(n-1)  \int   d(F_1(X_0,\ldots; \varepsilon_1),
 F_1(X_0,\ldots; y))P_{\varepsilon_1}(dy)  \\
 & \leq a_{ n-1}(n-1)  G_{ \varepsilon_1}(\varepsilon_1) \,,
\end{align*}
where the last line follows by \eqref{c2}. \hfill\qed

\vspace{0.3cm}

\noindent {\emph{Proof of Proposition  \ref{pr01}}.}
  By Proposition \ref{McD} and condition (\ref{BernsteinC}), it is  easy to see that, for any $k\in [1, n]$ and any $t \in [0, \delta^{-1})$,
\begin{eqnarray}\label{fin35}
\mathbb{E}\,[e^{t d_k } ] &=& 1+ \sum_{i=2}^{\infty} \frac{t^i}{i!}\, \mathbb{E}\,[   (d_k)^i ] \ \leq \ 1+ \sum_{i=2}^{\infty} \frac{t^i}{i!}\, \mathbb{E}\,[ | d_i|^i ] \nonumber\\
&\leq& 1+ \sum_{i=2}^{\infty} \frac{t^i}{i!}\, \big(a_{n-k}(n-k) \big)^i \mathbb{E}\,\Big[ \big(G_{\varepsilon_k}(\varepsilon_k) \big)^i \Big] \nonumber\\
&\leq& 1+ \sum_{i=2}^{\infty} \frac{t^i}{i!}\, \big(a_{n-k}(n-k) \big)^i \frac {i!}{2} V_k M^{i-2}
\leq 1+ \frac{t^2 V_k \big(   a_{n-k}(n-k)\big)^2 }{2 (1  -t\, \delta )}\, .
\end{eqnarray}
Using the inequality $1+t \leq e^t,$ we have, for any $k\in [1, n]$ and any $t \in [0, \delta^{-1})$,
\begin{eqnarray}
\mathbb{E}\,[e^{t d_k } ]
&\leq&  \exp \left (\frac{t^2V_k\big(a_{n-k}(n-k)\big)^2}{2 (1  -t\,\delta )}  \right ).
\end{eqnarray}
By the tower property of conditional expectation,  we deduce that, for any $k\in [1, n]$ and any $t \in [0, \delta^{-1})$,
\begin{eqnarray*}
\mathbb{E}\,\big[e^{ tS_n} \big]
&=&  \mathbb{E}\,\big[ e^ { tS_{ n-1}} \mathbb{E}\,  [e^ { td_n}  |\mathcal{F}_{n-1}  ] \big] \\
&\leq & \mathbb{E}\,\big[ e^ { tS_{ n-1}} \big] \exp \left( \frac{t^2V_n (a_0(0))^2 }{2 (1- t\,\delta)}  \right)\\
&\leq &  \exp \left(\frac{t^2 V}{2 (1- t\,\delta)} \right),
\end{eqnarray*}
which gives inequality (\ref{maindfs}).
Using the  Markov inequality, we obtain,
for any $x\geq 0$ and any $t \in [0, \delta^{-1})$,
\begin{eqnarray}
  \mathbb{P}\left( S_{n} \geq x \right)
 &\leq&  \mathbb{E}\, \big[e^{t\,(S_n -x) } \big ] \nonumber\\
 &\leq&   \exp \left(-t\,x  +  \frac{t^2 V}{2 (1- t\,\delta)}    \right)\, . \label{fines}
\end{eqnarray}
The minimum is reached at $$t=t(x):= \frac{2x/V}{ 2x\delta/V +1 + \sqrt{1+2x\delta/V} }\, .$$ Substituting $t=t(x)$ in (\ref{fines}),  we obtain the desired inequalities
\begin{eqnarray}
  \mathbb{P}\left( S_{n} \geq x \right)
 &\leq&  \exp \left( \frac{-x^2}{V(1+\sqrt{1+2x \delta/V})+x  \delta }\right)\, \nonumber\\
 &\leq&  \exp \left( \frac{-x^2}{2 (V +x  \delta ) }\right)\, ,\nonumber
\end{eqnarray}
where the last line follows from the inequality $\sqrt{1+2x\,\delta/V}  \leq 1+x\, \delta/V$. \hfill\qed

\vspace{0.3cm}

\noindent\emph{Proof of  Proposition   \ref{BerProp}.}
 By Taylor's expansion of $e^x$ and the fact that $\mathbb{E} S_{n}=0$, we have, for any $t\geq 0,$
 \begin{eqnarray}\label{sxsa}
 \mathbb{E}\bigg[\exp\bigg(t \frac{ S_{ n}}{\sqrt{n}} \bigg )\bigg] = 1+ \sum_{k=2}^{\infty} \frac{t^k}{k!}  \mathbb{E}\Big[\Big(\frac{ S_{ n}}{\sqrt{n}}\Big)^k \Big].
\end{eqnarray}
Using   Rio's inequality (see Theorem 2.1 of \cite{R09}):  we get, for any $k\geq 2,$
\begin{eqnarray}\label{R}
\Big(\mathbb{E}[|S_{n}|^k] \Big)^{2/k} \leq  (k-1)  \sum_{i=1}^n  \big(\mathbb{E}[|d_i|^k] \big)^{2/k}  ,
\end{eqnarray}
which is equivalent to
\begin{eqnarray}\label{ssxsfsd}
\mathbb{E}[|S_{n}|^k] \leq (k-1)^{k/2} \Big(\sum_{i=1}^n \big(\mathbb{E}[|d_i|^k] \big)^{2/k}\Big)^{k/2}.
\end{eqnarray}
Applying H\"{o}lder's  inequality  to  (\ref{ssxsfsd}), we deduce that, for any $k\geq 2,$
\begin{eqnarray}\label{sddd}
\mathbb{E}[|S_{n}|^k] \leq (k-1)^{k/2} n^{k/2-1} \sum_{i=1}^n  \mathbb{E}[|d_i|^k]  .
\end{eqnarray}
Applying the last inequality to (\ref{sxsa}), we have, for any $t\geq 0,$
 \begin{eqnarray}\label{bnlts}
\mathbb{E}\bigg[\exp\bigg(t \frac{ S_{ n}}{\sqrt{n}} \bigg )\bigg] \leq 1+ \sum_{k=2}^{\infty} \Big( \frac{t^k}{k!}  (k-1)^{k/2} n^{ -1 } \sum_{i=1}^n  \mathbb{E}[|d_i|^k] \Big).
\end{eqnarray}
 By point  2   of Proposition \ref{McD} and condition (\ref{Bercond02}), we deduce that,
for any integer $i\geq 1$,
\begin{eqnarray*}
\mathbb{E}[|d_i|^k] &\leq&  \mathbb{E}[|a_{n-i}(n-i)  H_{i,\varepsilon_i}(X_{i-1},...,X_1, \varepsilon_i)|^k ] \\
&\leq& \frac{1}{2} \frac{  k! \,( a_{n-1}(n-1) \epsilon )^{k-2}  }{(k-1)^{k/2}}   \mathbb{E}[ (a_{n-i}(n-i)  H_{i,\varepsilon_i}(X_{i-1},...,X_1, \varepsilon_i))^2],  \ \ \ \ \  k\geq  2.
\end{eqnarray*}
Hence  condition
(\ref{Bercond02}) implies that, for any $0\leq t < (a_{n-1}(n-1) \epsilon)^{-1}$,
 \begin{eqnarray}\label{ineq23}
\mathbb{E}\bigg[\exp\bigg(t \frac{ S_{ n}}{\sqrt{n}} \bigg )\bigg] \leq 1+ \sum_{k=2}^{\infty}  \frac{ \ \sigma_n^2}{2 } \, t^k (a_{n-1}(n-1) \epsilon )^{k-2} = 1+   \frac{t^2 \sigma_n^2  }{2\,(1- t a_{n-1}(n-1) \epsilon )} .
\end{eqnarray}
Using  the inequality $1+x \leq e^x  ,$ we have, for any $0\leq t < (a_{n-1}(n-1) \epsilon)^{-1}$,
 \begin{eqnarray*}
\mathbb{E}\bigg[\exp\bigg(t \frac{ S_{ n}}{\sqrt{n}} \bigg )\bigg] \ \leq \   \exp \Bigg( \frac{t^2  \sigma_n^2  }{2\,(1- t a_{n-1}(n-1)\epsilon)}   \Bigg).
\end{eqnarray*}
Applying Markov's inequality,   it is  easy to see that, for any $0\leq t < \sigma_n(a_{n-1}(n-1) \epsilon)^{-1}$ and any $ x \geq 0,$
\begin{eqnarray*}
  \mathbb{P}\left(   S_{ n} \geq x  V_n     \right)
   \leq   \exp\Big( - t x  \Big ) \mathbb{E}\left[   \exp\bigg(t \frac{ S_{ n}}{V_n } \bigg ) \right].
\end{eqnarray*}
Hence, it holds, for any $x>0,$
\begin{eqnarray*}
  \mathbb{P}\left(   S_{n} \geq x  V_n     \right)
    &\leq& \inf_{0\leq t <\sigma_n ( a_{n-1}(n-1) \epsilon)^{-1}}  \exp\bigg( - t x   + \frac{  t^2    }{2\,(1- t \, a_{n-1}(n-1)\epsilon/ \sigma_n)}    \bigg )\\
    &=& \exp\bigg( \frac{-\, x^2}{  1+ \sqrt{ 1+ 2x a_{n-1}(n-1) \epsilon
  /\sigma_n   }  +  x  a_{n-1}(n-1)\epsilon /\sigma_n   } \bigg),
\end{eqnarray*}
which gives   (\ref{jnsk01}).  Since
$   \sqrt{ 1+ 2x  a_{n-1}(n-1)\epsilon / \sigma_n    }  \leq  1+  x  a_{n-1}(n-1) \epsilon /  \sigma_n ,$
we get  (\ref{jnsk02})  from (\ref{jnsk01}).  \qed

\vspace{0.3cm}

\noindent\emph{Proof of Proposition  \ref{cram}.} Let $\delta=t_0/a_{n-1}(n-1).$  Since $\mathbb{E}  d_1  =0$, it is  easy to see that, for any $k \in [1, n]$ and any $t \in   [0, \delta )$,
\begin{eqnarray}\label{finsa13f}
\mathbb{E}\,[e^{t d_k } ] &=& 1+ \sum_{i=2}^{\infty} \frac{t^i}{i!}\, \mathbb{E}\,[   (d_k)^i ] \nonumber\\
&\leq& 1+ \sum_{i=2}^{\infty} \Big(\frac{t}{\delta} \Big)^i  \, \mathbb{E}\,\Big[\frac{1}{i!} | \delta  d_k |^i \Big].
\end{eqnarray}
Note that, for any $t\geq 0$,
\begin{align}\label{constantsecond}
\frac{t^i}{i !} e^{-t}  \leq \frac{i^i}{i !} e^{-i}
 \leq 2 e^{-2}, \quad \text{for any  $i\geq 2$,}
\end{align}
where the last line follows from the fact that $i^i e^{-i}/i!$ is decreasing in $i$. Note that the equality in \eref{constantsecond}
 is reached at $t=i=2$. Using \eref{constantsecond}, point 2 of Proposition \ref{McD} and condition (\ref{laplace}), we have, for any $i\geq2$ and any $k \in [1, n],$
\begin{eqnarray}\label{fisa3f}
\mathbb{E}\,\Big[\frac{1}{i!} |\delta d_k|^i \Big]   &=&  \Big(\frac {a_{n-k}(n-k)}{a_{n-1}(n-1)}\Big)^2 \mathbb{E}\,\Big[\frac{1}{i!} | \frac{t_0}{a_{n-k}(n-k)} d_k|^i \Big] \nonumber\\
 &\leq& \ 2e^{-2} \Big(\frac {a_{n-k}(n-k)}{a_{n-1}(n-1)}\Big)^2  \mathbb{E}\, [e^{ t_0|d_k| / a_{n-k}(n-k)}    ] \nonumber\\
 &\leq& \ 2e^{-2} \Big(\frac {a_{n-k}(n-k)}{a_{n-1}(n-1)}\Big)^2 {\mathbb E} \Big[ \exp \Big( t_0   G_{\varepsilon_k}(\varepsilon_k)\Big)\Big] \nonumber\\
 &\leq& \ 2e^{-2} \Big(\frac {a_{n-k}(n-k)}{a_{n-1}(n-1)}\Big)^2 K_k.
\end{eqnarray}
Combining the inequalities (\ref{finsa13f}) and (\ref{fisa3f}) together,  we obtain,  for any $t \in [0,   \delta )$,
\begin{eqnarray}\label{fin36}
  \mathbb{E}\,[e^{t d_k  }|\mathcal{F}_{k-1}] \  \leq \ \exp \left(\frac{2}{e^{2}} \frac{t^2 K_k \delta^{-2} }{1-t \delta^{-1}} \Big(\frac {a_{n-k}(n-k)}{a_{n-1}(n-1)}\Big)^2  \right) \, .
\end{eqnarray}
 By the tower property of conditional expectation, we deduce that, for any $t \in [0, \delta )$,
\begin{eqnarray}
\mathbb{E}\,\big[e^{ tS_n} \big] &=&  \mathbb{E}\,\big[ \,\mathbb{E}\, [e^{ t S_n} |\mathcal{F}_{n-1}  ] \big]  =   \mathbb{E}\,\big[ e^ { t S_{n-1}} \mathbb{E}\,  [e^ { t d_n}  |\mathcal{F}_{n-1}  ] \big]\nonumber\\
\nonumber \\
&\leq & \mathbb{E}\,\big[  e^ { t S_{n-1}}  \big] \exp \left( \frac{2}{e^{2}} \frac{t^2 K_k \delta^{-2} }{1-t \delta^{-1}} \Big(\frac {a_{0}(0)}{a_{n-1}(n-1)}\Big)^2 \right)  \nonumber\\
&\leq &  \exp \left( \frac{t^2 K  \delta^{-2} }{1-t \delta^{-1}}   \right),
\end{eqnarray}
where $$K=\frac{2}{e^{2}}  \sum_{k=1}^{n} \Big(\frac {a_{n-k}(n-k)}{a_{n-1}(n-1)}\Big)^2 K_k  .$$
Using the exponential Markov  inequality, it follows that,
for any $x\geq 0$ and any $t \in [0, \delta )$,
\begin{eqnarray}
  \mathbb{P}\left( S_{n} \geq x \right)
 &\leq&  \mathbb{E}\, [e^{t\,(S_n -x) }  ] \nonumber\\
 &\leq&   \exp \left(-t x  +  \frac{t^2 K  \delta^{-2} }{1-t \delta^{-1}}   \right)\, .  \label{fisaxs}
\end{eqnarray}
The  minimum is reached at $$t=t(x):= \frac{x\delta^2/K}{x\delta/K + 1 + \sqrt{1+x\delta/K}} .$$ Substituting $t=t(x)$ in (\ref{fisaxs}),  we obtain the desired inequality (\ref{Berie3}).  The second desired inequality follows by the fact $\sqrt{1+ x \delta / K } \leq  1+  x \delta / 2K   $.    \hfill\qed

\vspace{0.3cm}

\noindent\emph{Proof of Proposition \ref{findsa}.}  Set $\xi_k= d_k/a_{n-1}(n-1)$ for any $ k \in [1, n].$ Denote $\xi^+=\max{\{\xi, 0 \}}$.  Using  point 2 of Proposition \ref{McD} and condition (\ref{laplace2}), we have, for any $k \in [1, n]$,
\begin{eqnarray}\label{finsa3f}
\sum_{k=1}^{n}\mathbb{E}\,[ \xi_k^2 e^{ (\xi_k^+)^p } | \mathcal{F}_{k-1} ] &=& \sum_{k=1}^{n} \Big(\frac { 1}{a_{n-1}(n-1)}\Big)^2 \mathbb{E}\,[ d_k^2 e^{ (d_k^+/a_{n-1}(n-1) )^p } | \mathcal{F}_{k-1} ]\nonumber\\
&\leq& \sum_{k=1}^{n} \Big(\frac { a_{n-k}(n-k)}{a_{n-1}(n-1)}\Big)^2 \mathbb{E}\,\Big[ (G_{\varepsilon_k}(\varepsilon_k))^2 \exp\Big ( \big(\frac{a_{n-k}(n-k)}{a_{n-1}(n-1)} G_{\varepsilon_k}(\varepsilon_k)  \big)^p \Big ) \Big| \mathcal{F}_{k-1} \Big]\nonumber\\
&\leq& \sum_{k=1}^{n} \Big(\frac { a_{n-k}(n-k)}{a_{n-1}(n-1)}\Big)^2 \mathbb{E}\,[ (G_{\varepsilon_k}(\varepsilon_k))^2 \exp\big ( (  G_{\varepsilon_k}(\varepsilon_k)  )^p \big ) | \mathcal{F}_{k-1} ]\nonumber\\
&\leq& \sum_{k=1}^{n} \Big(\frac { a_{n-k}(n-k)}{a_{n-1}(n-1)}\Big)^2  K_k=: K.
\end{eqnarray}
 Using Theorem 2.1 of  \cite{FGL17}, we obtain,  for any $x > 0$,
\begin{eqnarray}
 &&  \mathbb{P}\Big(  S_n \geq x \, a_{n-1}(n-1)  \Big) \ \ \ \ \  \ \ \  \nonumber \\
&&  \leq \left\{ \begin{array}{ll}
\exp\Bigg(\displaystyle -  \frac{x^2}{2K }   \Bigg) +  K \bigg( \frac{ x  } {K }\bigg)^{2/(1-p)}\exp\Bigg( -\Big(\frac{K}{x  } \Big)^{p/(1-p)} \Bigg)    & \textrm{if\ \  $0\leq x < K^{1/(2-p)}$}   \\
\vspace{-0.2cm}\\
\exp \Bigg(\displaystyle - x^p\Big( 1- \frac{K}{2\, x^{2-p}}\Big)  \Bigg) + K\frac{1}{x^2}\exp \bigg(- x^{p}   \bigg)   & \textrm{if\ \  $x \geq K^{1/(2-p)}$ ,}
\end{array} \right.
\end{eqnarray}
and moreover, if $K\geq 1,$ then for any $x>0,$
\begin{eqnarray}
  \mathbb{P}\Big(  S_n \geq x \, a_{n-1}(n-1)  \Big)   \leq  2\exp \Bigg(- \frac{x^2 }{2(K+ x^{2-p})} \Bigg) .
\end{eqnarray}
The last  inequity is equivalent to our desired inequality.     \hfill\qed

\vspace{0.3cm}

\noindent\emph{Proof of Proposition \ref{propExp}.} Clearly,  by point 2 of Proposition \ref{McD} and condition (\ref{Bernsteinmoment}), it holds:
   for any  $k \in [1, n],$
\begin{eqnarray}\label{fin36}
\mathbb{E}\,\left[ \exp\left\{  | a_{n-1}^{-1}(n-1) d_k|^{\frac{2\alpha}{1-\alpha}} \right\} \right]  \leq   C_1 .
\end{eqnarray}
Applying Theorem 2.1 of  \cite{Fx1} to martingale sequence
 $( a_{n-1}^{-1} (n-1)d_k, \mathcal{F}_k)_{k=1,..,n}$,   we obtain  the desired inequality.
   \hfill\qed

\vspace{0.3cm}

\noindent\emph{Proof of Proposition \ref{po71}.} We apply Corollary 2.3 of  \cite{FGL12} with the truncature level $y a_{n-1}(n-1)$.
By point 3 of Proposition \ref{McD},   $|d_k| \leq a_{n-k}(n-k) G_{\varepsilon_k}(\varepsilon_k)$ for any $k\in  [1, n]$. Hence, for any $k\in [1, n]$,
$$
 {\mathbb E}\big[d_k^2 {\bf 1}_{\{d_k \leq y a_{n-k}(n-k)\}}|{\mathcal F}_{i-1}\big]
  \leq \Big( a_{n-k}(n-k)\Big)^2
 {\mathbb E} \big[  \big( G_{\varepsilon_k}(\varepsilon_k)\big)^2\big] \leq
  \Big(a_{n-k}(n-k)\Big)^2 V_k\, .
$$
Then  it follows from  Corollary 2.3 of   \cite{FGL12} that
$$
{\mathbb P}(S_n>x)  \leq    H_n \left(\frac{x }{y a_{n-1}(n-1)} , \frac{\sqrt{V} }{y a_{n-1}(n-1) }\right) \\
+ {\mathbb P}\bigg( \max_{1 \leq k \leq n} d_k >  y a_{n-1}(n-1) \bigg)\, .
$$
Inequality \eref{Hineq} follows by applying point 3  of
Proposition \ref{McD} again. \hfill \qed

\vspace{0.3cm}

\noindent\emph{Proof of Proposition \ref{RCfG}.}
By Proposition \ref{McD}  and condition (\ref{sdvcx}), it follows that
\begin{eqnarray*}
\sum_{k=1}^n\mathbb{E}[|d_k|^p | \mathcal{F}_{k-1} ] &\leq&   \sum_{k=1}^n\mathbb{E}[| a_{n-k}(n-k) G_{\varepsilon_k}(\varepsilon_k)|^p ] \\
&\leq &  \sum_{k=1}^n \big( a_{n-k}(n-k) \big)^pA_k(p)=A(p).
\end{eqnarray*}
Notice that $A(2)=V$. Using Corollary $3'$ of   \cite{F73}, we obtain the desired inequality. \hfill\qed

\vspace{0.3cm}

\noindent\emph{Proof of Proposition \ref{classic}.}
Let
$$
u_{k-1}(x_1, \ldots, x_{k-1})= \text{ess\,inf}_{(x_{n-i})_{i\geq n}, \varepsilon_k} g_k(x_1, \ldots, x_{k-1}, F_k((x_{n-i})_{i\geq 1};  \varepsilon_k))
$$
and
$$
v_{k-1}(x_1, \ldots, x_{k-1})= \text{ess\,sup}_{(x_{n-i})_{i\geq n}, \varepsilon_k}  g_k(x_1, \ldots, x_{k-1}, F_k((x_{n-i})_{i\geq 1};  \varepsilon_k))
$$
From the proof of Proposition \ref{McD},  it follows that
$$
u_{k-1}(X_1, \ldots, X_{k-1}) \leq d_k \leq v_{k-1}(X_1, \ldots, X_{k-1})\, .
$$
By Proposition \ref{McD} and  condition (\ref{in850}), we have
\[
v_{k-1}(X_1, \ldots, X_{k-1})-u_{k-1}(X_1, \ldots, X_{k-1}) \leq a_{n-k}(n-k) M_k\, .
\]
Now, following exactly the proof of  Theorem 3.1 of   \cite{R13} with $\Delta_k =  a_{n-k}(n-k) M_k$
we  obtain the inequalities (\ref{rio1}) and (\ref{rio2}).  Since
for any $x \in [0,1)$,
$\ell^*(x) \geq (x^2-2x) \ln(1-x)$, inequality (\ref{riosbound})
follows from (\ref{rio2}). \hfill\qed

\vspace{0.3cm}

\noindent\emph{Proof of Proposition \ref{VBEI}.} Using an improvement of the von Bahr-Esseen inequality  (cf.\ inequality (1.11) in   \cite{P10}), we get
$$
\| S_n \|_p^p \leq \|d_1\|_p^p +  \tilde C_p \sum_{k=2}^n  \|\, d_k \|_p^p \, ,
$$
where $\tilde C_p$ is a constant satisfying $\tilde C_p \leq 2^{2-p}$
for any $p \in [1,2]$, and it is described in Proposition 1.8 of \cite{P10}.
By Proposition \ref{McD}, we have
\begin{eqnarray*}
\| S_n \|_p^p   &\leq&  \, \bigg( \,  \big(a_{n-1}(n-1)\big)^p   {\mathbb E} \Big[  \big( H_{1, \varepsilon_1}(  \varepsilon_1)\big)^p\Big]    + \tilde C_p
    \sum_{k=2}^n \big(a_{n-k}(n-k)\big)^p {\mathbb E} \Big[  \big( H_{k, \varepsilon_k}(X_{k-1},...,X_1, \varepsilon_k)\big)^p\Big]\, \bigg) \\
    &\leq&  \bigg(\big(a_{n-1}(n-1)\big)^p A_1(p)   +
\tilde C_p  \sum_{k=2}^n \big(a_{n-k}(n-k)\big)^p A_k(p) \, \bigg),
\end{eqnarray*}
which gives the desired inequality. \hfill\qed

\vspace{0.3cm}

\noindent\emph{Proof of Proposition \ref{weakVBEI}.} By Proposition 3.3 of   \cite{CDM17}, we deduce that, for  any $x>0,$
\begin{equation}\label{inew01}
 {\mathbb P}(|S_{n}|\geq x ) \leq  \frac{  C_p}{x^p}\, \sum_{k=1}^n   \left\| d_k\right\|_{w,p}^p  .
\end{equation}
From point 2 of Proposition \ref{McD} and condition (\ref{weakVB}),  it follows that,
for any  $k \in [1, n],$
\begin{equation}\label{inew03}
 \left \| d_k\right\|_{w,p}^p  \leq   \left \|  a_{n-k}(n-k) H_{k,\varepsilon_k}(X_{k-1},...,X_1, \varepsilon_k)    \right \|_{w,p}^p  \leq   \left( a_{n-k}(n-k) \right)^p   A_k(p).
\end{equation}
Combining  \eref{inew01}  and  \eref{inew03} together, we get the desired inequality.  \hfill\qed

\end{document}